\documentclass[11pt, a4paper]{article}
\usepackage{amsmath}
\usepackage{amsfonts}
\usepackage{amssymb}
\usepackage{ulem}

\usepackage{amsthm}
\usepackage{url}
\usepackage{dcolumn}

\usepackage{cite}
\usepackage{fancyhdr}

\usepackage{graphicx, color}
\usepackage{hyperref}
\usepackage{xcolor}

\newtheorem*{rem*}{Remark}

\begin{document}

\title{Convex pentagons and convex hexagons that can form \\rotationally 
symmetric tilings}
\author{ Teruhisa SUGIMOTO$^{ 1), 2)}$ }
\date{}
\maketitle

{\footnotesize

\begin{center}
$^{1)}$ The Interdisciplinary Institute of Science, Technology and Art

$^{2)}$ Japan Tessellation Design Association

E-mail: ismsugi@gmail.com
\end{center}

}

{\small
\begin{abstract}
\noindent
In this study, the properties of convex hexagons that can form 
rotationally symmetric edge-to-edge tilings are discussed. Because the 
convex hexagons are equilateral convex parallelohexagons, convex 
pentagons generated by bisecting the hexagons can form rotationally 
symmetric non-edge-to-edge tilings. In addition, under certain 
circumstances, tiling-like patterns with an equilateral convex polygonal 
hole at the center can be formed using these convex hexagons or pentagons.
\end{abstract}
}


%%%%%%%%%%%%%%%%%%%%%%%%%%%%%%%%%%%%%%%%%%%%%%%%%%%%%%%%%%%%%%%%%%%%%%
%%%%%%%%%%%%%%%%%%%%%%%%%%%%%%%%%%%%%%%%%%%%%%%%%%%%%%%%%%%%%%%%%%%%%%

\section{Introduction}
\label{section1}

In~\cite{Klaassen_2016}, Klaassen has demonstrated that there exist countless 
rotationally symmetric non-edge-to-edge tilings with convex pentagonal 
tiles.\footnote{ 
A \textit{tiling} (or \textit{tessellation}) of the plane is a collection of sets 
that are called tiles, which covers a plane without gaps and overlaps, except 
for the boundaries of the tiles. The term ``tile" refers to a topological disk, 
whose boundary is a simple closed curve. If all the tiles in a tiling are of the 
same size and shape, then the tiling is \textit{monohedral}~\cite{G_and_S_1987, 
wiki_pentagon_tiling}. In this study, a polygon that admits a monohedral tiling 
is called a \textit{polygonal tile}~\cite{Sugimoto_NoteTP}. Note that, in monohedral 
tiling, it admits the use of reflected tiles.
}$^{,}$\footnote{ 
A tiling by convex polygons is \textit{edge-to-edge} if any two convex polygons 
in a tiling are either disjoint or share one vertex or an entire edge in common. 
Then other case is \textit{non-edge-to-edge}~\cite{G_and_S_1987, Sugimoto_NoteTP}.
} The convex pentagonal tiles are considered to be equivalent to bisecting 
equilateral convex parallelohexagons, which are hexagons where the opposite 
edges are parallel and equal in length. Thus, there exist countless rotationally 
symmetric edge-to-edge tilings with convex hexagonal tiles.

\renewcommand{\figurename}{{\small Figure.}}
\begin{figure}[htb]
 \centering\includegraphics[width=13.5cm,clip]{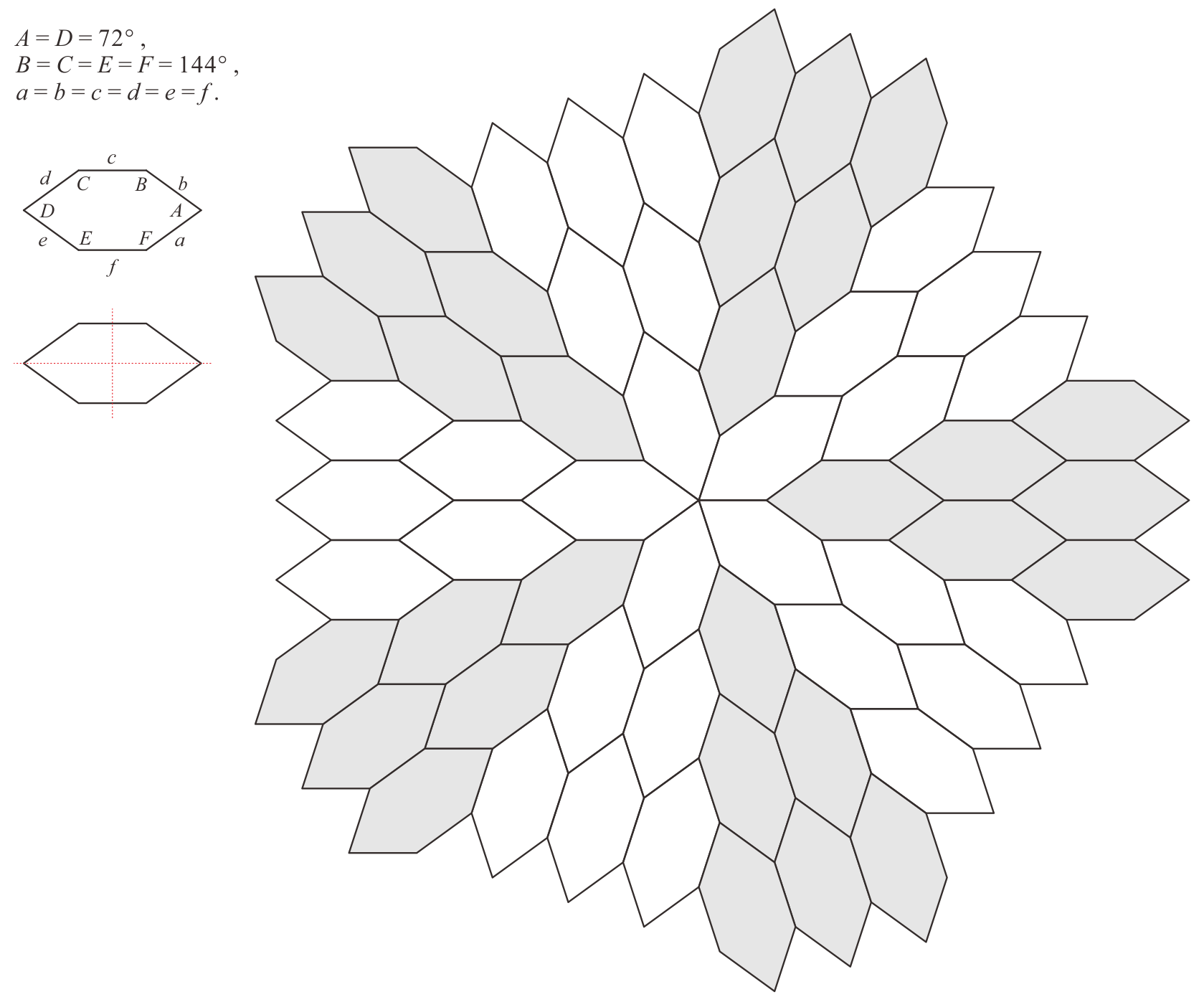} 
  \caption{{\small 
Rotationally symmetric tiling with $D_{5}$ symmetry that formed by 
a convex hexagonal tile with $D_{2}$ symmetry (Note that the gray area in the 
figure is used to clearly depict the structure)
} 
\label{fig01}
}
\end{figure}

\begin{figure}[!h]
 \centering\includegraphics[width=15cm,clip]{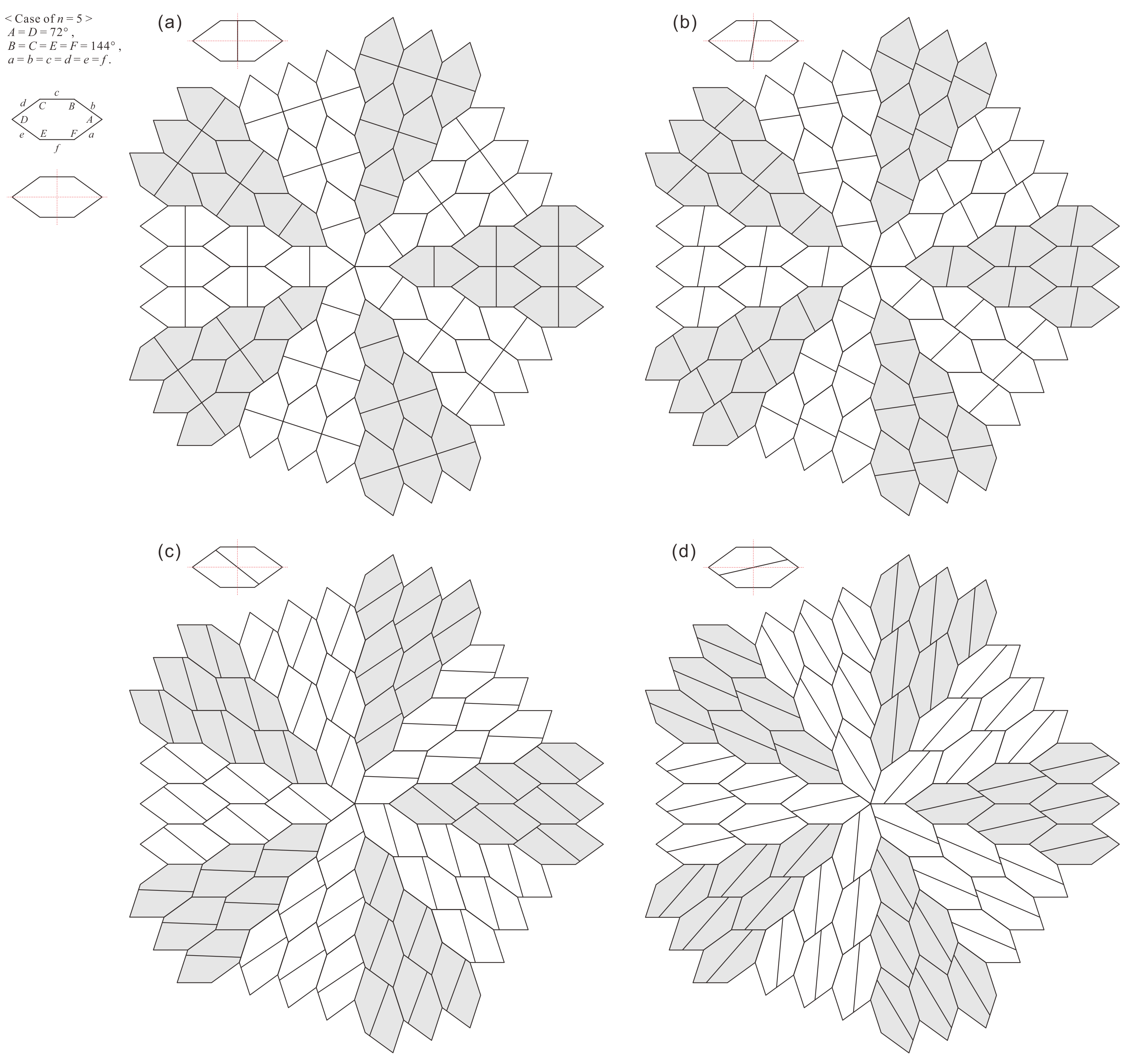} 
  \caption{{\small 
Five-fold rotationally symmetric non-edge-to-edge tilings 
with convex pentagons
} 
\label{fig02}
}
\end{figure}

\begin{figure}[htbp]
 \centering\includegraphics[width=15cm,clip]{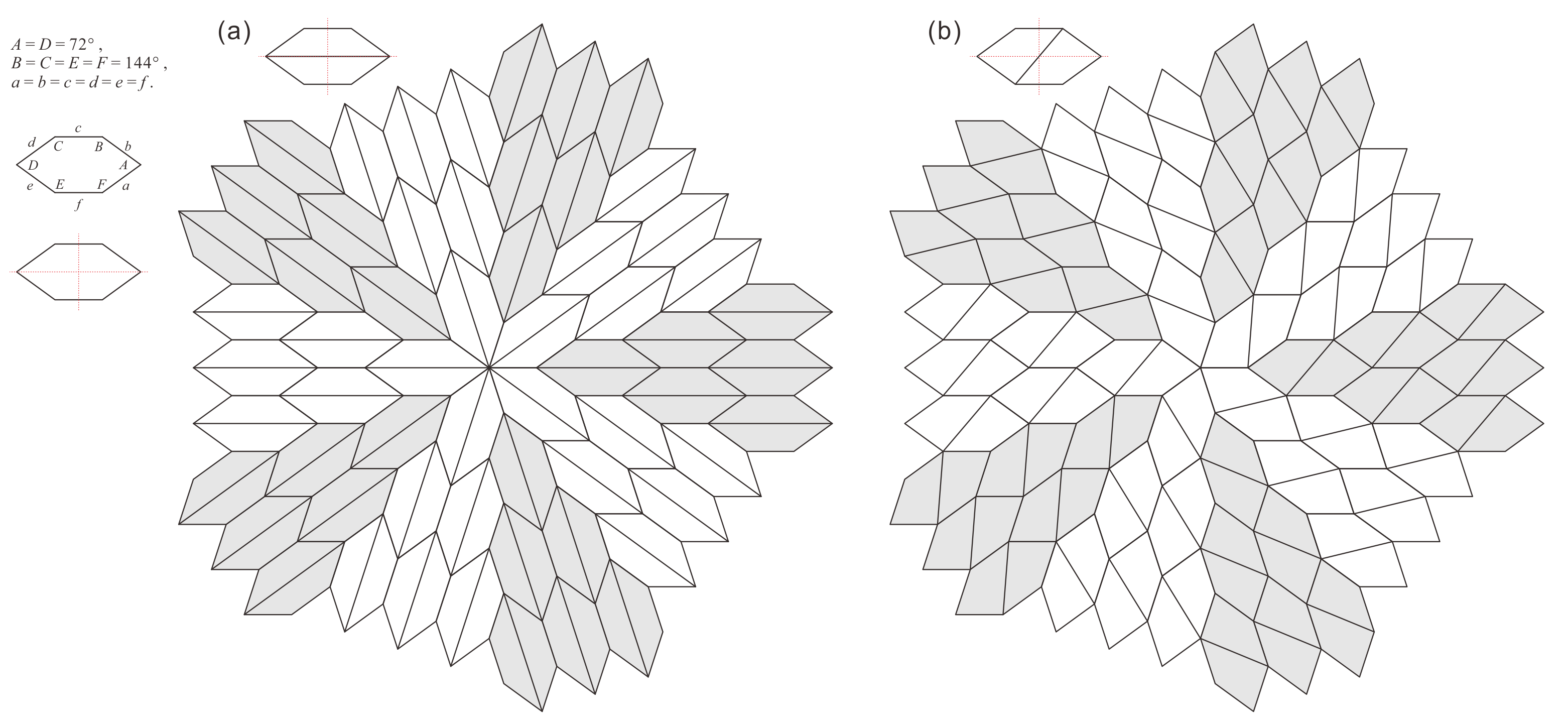} 
  \caption{{\small 
Five-fold rotationally symmetric edge-to-edge tilings 
with convex quadrangles
} 
\label{fig03}
}
\end{figure}

Figure~\ref{fig01} shows a five-fold rotationally symmetric edge-to-edge tiling by a 
convex hexagonal tile (equilateral convex parallelohexagon) that satisfies the conditions, 
``$A = D = 72^ \circ ,\; B = C = E = F = 144^ \circ ,\;a = b = c = d = e = f$." 
Figure~\ref{fig02} shows examples of five-fold rotationally symmetric non-edge-to-edge 
tilings with convex pentagonal tiles generated by bisecting the equilateral convex 
parallelohexagon in Figure~\ref{fig01}. Because the equilateral convex parallelohexagons 
have two-fold rotational symmetry, the number of ways to form convex pentagons 
generated by bisecting the convex hexagons (the dividing line needs to be a 
straight line that passes through the rotational center of the convex hexagon 
and intersects the opposite edge) is countless. The bisecting method can be 
divided into three cases: 
Case (i) the dividing line intersects edges $c$ and $f$, as shown in 
Figures~\ref{fig02}(a) and \ref{fig02}(b); 
Case (ii) the dividing line intersects edges $a$ and $d$, as shown in Figure~\ref{fig02}(c); 
Case (iii) the dividing line intersects the edges $b$ and $e$, as shown in 
Figure~\ref{fig02}(d). 
If the dividing line is selected such that the opposite vertices of the equilateral 
convex parallelohexagons are connected, the parallelohexagons contain two 
congruent convex quadrangles. In such cases, as shown in Figure~\ref{fig03}, 
there exist five-fold rotationally symmetric edge-to-edge tilings with 
convex quadrangles.

\renewcommand{\figurename}{{\small Figure.}}

\begin{figure}[htbp]
 \centering\includegraphics[width=15cm,clip]{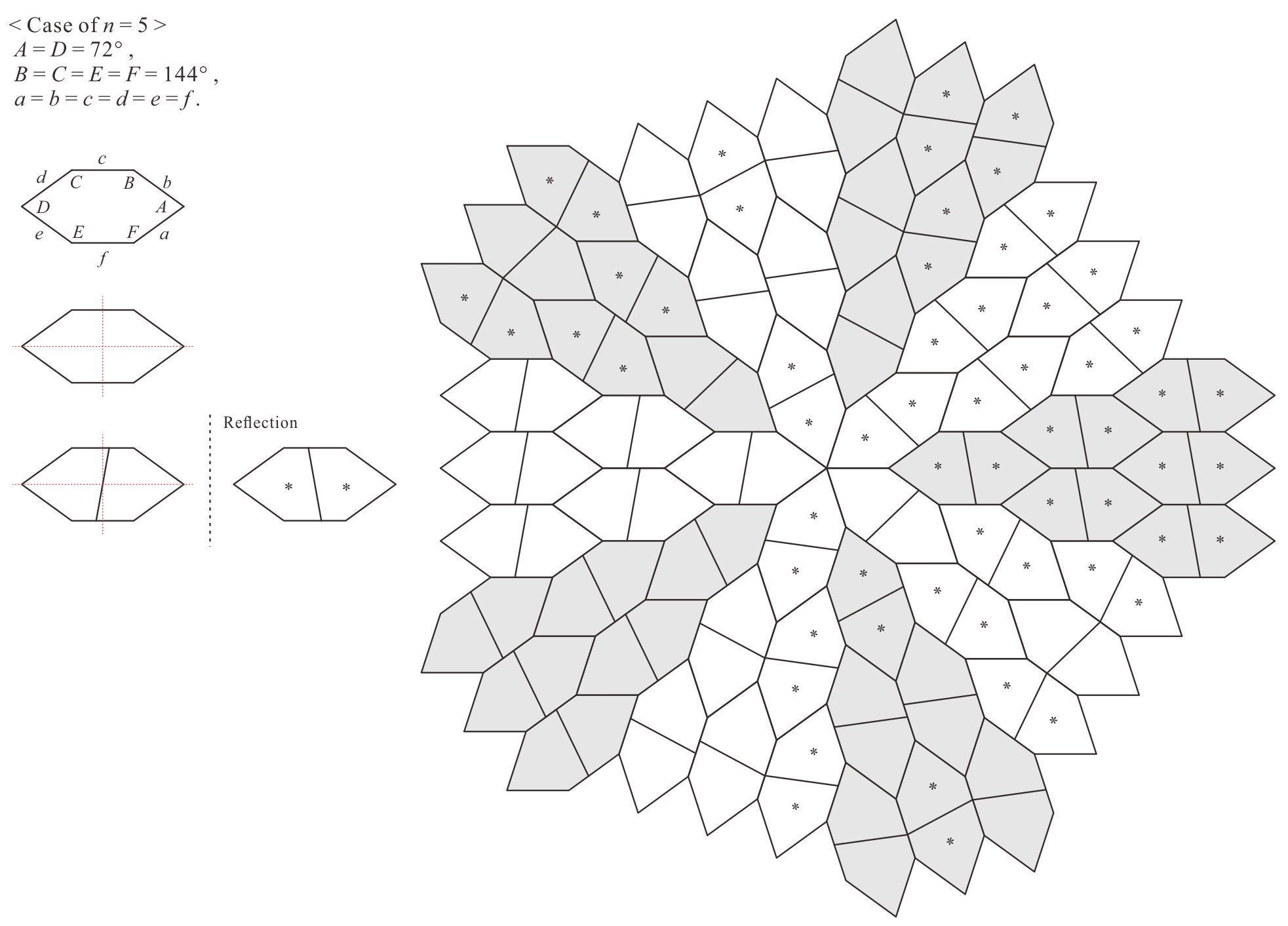} 
  \caption{{\small 
Random tiling based on the five-fold rotationally symmetric 
tiling structure of a convex pentagonal tile
} 
\label{fig04}
}
\end{figure}

The equilateral convex parallelohexagons that satisfy ``$A = D \ne B = C =  E 
= F,\;a = b = c = d = e = f$," as demonstrated in Figure~\ref{fig01}, have two-fold 
rotational symmetry and two axes of reflection symmetry passing through the 
center of the rotational symmetry (hereafter, this property is described as 
$D_{2}$ symmetry\footnote{ 
``$D_{2}$" is based on the Schoenflies notation for symmetry in a two-dimensional 
point group~\cite{wiki_point_group, wiki_schoenflies_notation}. ``$D_{n}$" represents 
an $n$-fold rotation axis with $n$ reflection symmetry axes. The notation for 
symmetry is based on that presented in \cite{Klaassen_2016}.
}). Therefore, the parallelohexagon 
and reflected parallelohexagon have identical outlines. (In the parallelohexagon 
with $D_{2}$ symmetry, Cases (ii) and (iii) can be regarded as having a 
reversible relationship.) If two convex pentagons are generated in the 
equilateral convex parallelohexagon with $D_{2}$ symmetry, which is bisected 
by a dividing line that does not overlap with the axis of reflection symmetry, 
as shown in Figures~\ref{fig02}(b), \ref{fig02}(c), and \ref{fig02}(d), the reflected 
parallelohexagon has the same outline. However, the arrangement of the inner 
convex pentagons is different. By using this property, the reflected convex 
pentagons can be freely incorporated into the tiling. Figure~\ref{fig04} shows a 
random tiling based on the five-fold rotationally symmetric tiling structure 
of a convex pentagonal tile created using this property.

In this study, the properties of convex hexagons and pentagons that can 
form rotationally symmetrical tilings presented by Klaassen are explored. 
The convex hexagons and pentagons can form rotationally symmetric 
tilings and rotationally symmetric tiling-like patterns with an equilateral 
convex polygonal hole at the center. Note that the tiling-like patterns are not 
considered tilings due to the presence of a gap, but are simply called 
tilings in this study. Herein, the various types of convex hexagons and 
pentagons are introduced and explored.

\renewcommand{\figurename}{{\small Figure.}}

\renewcommand{\figurename}{{\small Figure.}}

%%%%%%%%%%%%%%%%%%%%%%%%%%%%%%%%%%%%%%%%%%%%%%%%%%%%%%%%%%%%%%%%%%%%%%
%%%%%%%%%%%%%%%%%%%%%%%%%%%%%%%%%%%%%%%%%%%%%%%%%%%%%%%%%%%%%%%%%%%%%%

\section{Rotationally symmetric tilings}
\label{section2}

In~\cite{Klaassen_2016}, the theorem and proof that countless non-edge-to-edge 
tilings can be formed with convex pentagonal tiles are presented, and the 
figures of five-fold and seven-fold rotationally symmetric tilings are shown. 
The conditions of the convex hexagonal tiles that can be formed with an 
$n$-fold rotationally symmetric edge-to-edge tiling are expressed in 
(\ref{eq1}), while considering the pentagonal tilings as convex hexagonal tilings. 
Note that the vertices (interior angles) and edges of the convex hexagon will 
be referred to using the nomenclature shown in Figure~\ref{fig01}.

\begin{equation}
\label{eq1}
\left\{ {\begin{array}{l}
 A = D = \dfrac{360^ \circ }{n\strut}, \\ 
 B = E, \\ 
 C = F, \\ 
 A + B + C = 360^ \circ , \\ 
 a = b = c = d = e = f, \\ 
 \end{array}} \right.
\end{equation}

\noindent
where $n$ is an integer greater than or equal to three, because 
$0^ \circ < A < 180^ \circ $. In (\ref{eq1}), the $A$ and $D$ pair are 
selected as the vertices for $\frac{360^ \circ }{n}$. 
However, due to the symmetry of the equilateral convex parallelohexagon, the 
$B$ and $E$ pair, as well as the $C$ and $F$ pair, are also possible 
vertex symbol starting points.

When convex hexagonal tiles satisfy (\ref{eq1}), where $B = C = E = F$, they 
become equilateral convex parallelohexagons with $D_{2}$ symmetry, as shown 
in Figure~\ref{fig01}. The conditions of the convex hexagonal tiles, in this case, are 
expressed in (\ref{eq2}).

\begin{equation}
\label{eq2}
\left\{ {\begin{array}{l}
 A = D = \dfrac{360^ \circ }{n\strut}, \\ 
 B = C = E = F = 180^ \circ -\dfrac{A}{2} = 180^ \circ -\dfrac{180^ \circ}{n\strut}, \\ 
 a = b = c = d = e = f. \\ 
 \end{array}} \right.
\end{equation}

The convex hexagons that satisfy (\ref{eq1}) belong to the Type 1 
family.\footnote{
It is known that convex hexagonal tiles belong to at least 
one of the three families referred to as a ``Type"~\cite{G_and_S_1987, 
Wolf_hexagon_tiling}.
} In contrast, the convex hexagons that satisfy (\ref{eq2}) belong to the 
Type 1 and Type 2 families. When $n = 3$ in (\ref{eq2}), a regular hexagon has 
$A = B = C = D = E = F = 120^ \circ $, which belongs to the Type 1, Type 2, 
and Type 3 families. The convex pentagons generated by bisecting 
the convex hexagons belong to the Type 1 family.\footnote{
To date, fifteen families of convex pentagonal tiles, each of them 
referred to as a ``Type," are known~\cite{G_and_S_1987, Sugimoto_NoteTP, 
wiki_pentagon_tiling}. For example, if the sum of three consecutive angles 
in a convex pentagonal tile is  $360^ \circ $, the pentagonal tile belongs 
to the Type 1 family. Convex pentagonal tiles belonging to some families 
also exist. Known convex pentagonal tiles can form periodic tiling. 
In May 2017, Micha\"{e}l Rao declared that the complete list of Types 
of convex pentagonal tiles had been obtained (i.e., they have only the 
known 15 families), but it does not seem to be fixed as of March 
2020~\cite{wiki_pentagon_tiling}.
}

\begin{table}[htb]
 \begin{center}
{\small
\caption[Table 1]{
Interior angles of convex hexagons satisfying (\ref{eq2}) that can 
form the $n$-fold rotationally symmetric edge-to-edge tilings
}
\label{tab1}
}
\begin{tabular}
{c| D{.}{.}{2} D{.}{.}{2} D{.}{.}{2} D{.}{.}{2} D{.}{.}{2} D{.}{.}{2} |c}
%{c|rrrrrr|c}
\hline
\raisebox{-1.50ex}[0cm][0cm]{$n$}& 
\multicolumn{6}{c|}{\shortstack{ Value of interior angle (degree) } } & 
\raisebox{-2.6ex}[0.7cm][0.5cm]{\shortstack{Figure \\number}} \\

 & 
\textit{A} & 
\textit{B}& 
\textit{C}& 
\textit{D}& 
\textit{E}& 
\textit{F}& 
  \\
\hline
3& 
120& 
120& 
120& 
120& 
120& 
120& 
\ref{fig15}\\
\hline
4& 
90& 
135& 
135& 
90& 
135& 
135& 
\ref{fig16}\\
\hline
5& 
72& 
144& 
144& 
72& 
144& 
144& 
\ref{fig01}, \ref{fig02}\\
\hline
6& 
60& 
150& 
150& 
60& 
150& 
150& 
\ref{fig17}\\
\hline
7& 
51.43 & 
154.29 & 
154.29 & 
51.43 & 
154.29 & 
154.29 & 
\ref{fig18}\\
\hline
8& 
45& 
157.5& 
157.5& 
45& 
157.5& 
157.5& 
\ref{fig19}\\
\hline
9& 
40& 
160& 
160& 
40& 
160& 
160& 
 \\
\hline
10& 
36& 
162& 
162& 
36& 
162& 
162& 
 \\
\hline
11& 
32.73 & 
163.64 & 
163.64 & 
32.73 & 
163.64 & 
163.64 & 
 \\
\hline
12& 
30& 
165& 
165& 
30& 
165& 
165& 
 \\
\hline
13& 
27.69 & 
166.15 & 
166.15 & 
27.69 & 
166.15 & 
166.15 & 
 \\
\hline
14& 
25.71 & 
167.14 & 
167.14 & 
25.71 & 
167.14 & 
167.14 & 
 \\
\hline
15& 
24& 
168& 
168& 
24& 
168& 
168& 
 \\
\hline
16& 
22.5& 
168.75& 
168.75& 
22.5& 
168.75& 
168.75& 
 \\
\hline
17& 
21.18 & 
169.41 & 
169.41 & 
21.18 & 
169.41 & 
169.41 & 
 \\
\hline
18& 
20& 
170& 
170& 
20& 
170& 
170& 
 \\
\hline
...& 
...& 
...& 
...& 
...& 
...& 
...& 
 \\
\hline
\end{tabular}

\end{center}

\end{table}

\renewcommand{\figurename}{{\small Figure.}}
\begin{figure}[!h]
 \centering\includegraphics[width=15cm,clip]{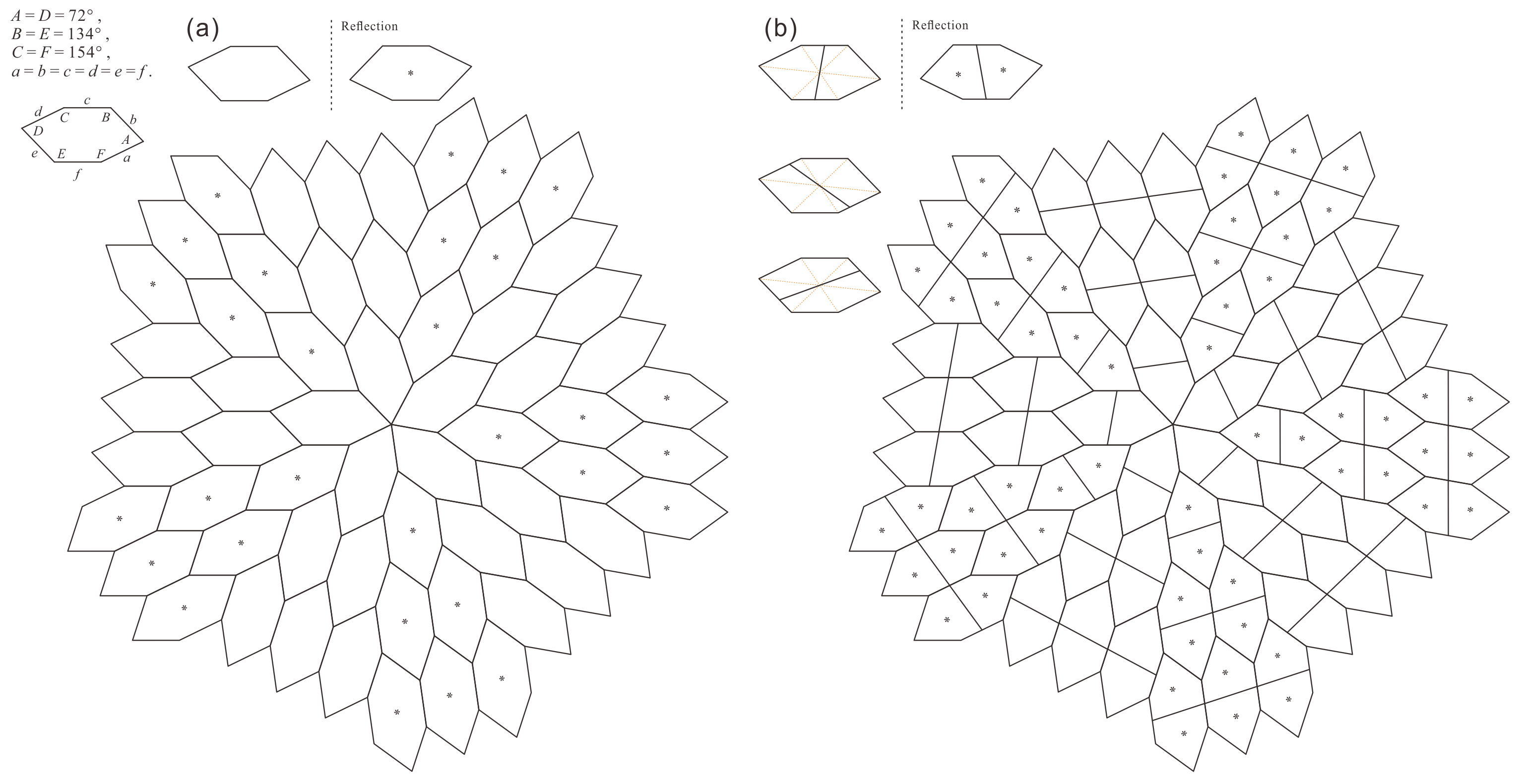} 
  \caption{{\small 
Rotationally symmetric edge-to-edge tiling with  $C_{5}$ symmetry 
by a convex hexagonal tile that satisfies ``$A = D = 72^ \circ,\; 
B = E = 134^ \circ ,\;C = F = 154^ \circ ,\; a = b = c = d = e = f$," 
and the five-fold rotationally symmetric non-edge-to-edge tiling by 
a convex pentagonal tile generated by bisecting the convex hexagon
} 
\label{fig05}
}
\end{figure}

Table~\ref{tab1} presents some of the relationships between the interior angles of 
convex hexagons satisfying (\ref{eq2}) that can form the $n$-fold rotationally 
symmetric edge-to-edge tilings. (For $n = 3\!-\!8$, tilings with convex 
hexagonal tiles and convex pentagonal tiles generated by bisecting them are 
drawn. For further details, Figures~\ref{fig01}, \ref{fig02}, \ref{fig15}--\ref{fig19}.) 
As mentioned above, when $n = 3$ in (\ref{eq2}), the convex hexagonal tile becomes a regular 
hexagon. Owing to the $D_{6}$ symmetry of the regular hexagon, the pentagonal 
pair corresponding to the bisected regular hexagon can be arranged by freely 
combining operations of $120^ \circ $ or $240^ \circ $ turnings and their reflections. 
(In \cite{Iliev_2018}, there are figures that depict a few tilings with such 
operations.)

For $n = 3$, the convex hexagonal tile that satisfies (\ref{eq2}) is a regular 
hexagon, and thus, its tiling has $D_{6}$ symmetry. For $n \ge 4$ each 
$n$-fold, rotationally symmetric edge-to-edge tiling by a convex hexagonal 
tile that satisfies (\ref{eq2}) has $D_{n}$ symmetry~\cite{Klaassen_2016}. 
Therefore, the tiling in Figure~\ref{fig01} has $D_{5}$ symmetry. 
The $n$-fold rotationally symmetric tilings with convex pentagons or convex 
quadrangles generated by bisecting the convex hexagons that satisfy 
(\ref{eq2}) along the axis of reflection symmetry have $D_{n}$ symmetry. 
Therefore, the tilings of Figures~\ref{fig02}(a) and \ref{fig03}(a) have 
$D_{5}$ symmetry.

When the convex hexagonal tile that satisfies (\ref{eq1}) has $B = E \ne C = F$, 
the equilateral convex parallelohexagon has two-fold rotational symmetry 
but no axis of reflection symmetry (hereafter, this property is described as 
$C_{2}$ symmetry\footnote{
``$C_{2}$" is based on the Schoenflies notation for symmetry in a two-dimensional 
point group~\cite{wiki_point_group, wiki_schoenflies_notation}. ``$C_{n}$" represents 
an $n$-fold rotation axis without reflection.
}). The $n$-fold rotationally symmetric edge-to-edge 
tilings by equilateral convex parallelohexagons with $C_{2}$ symmetry have $C_{n}$ 
symmetry because they have rotational symmetry but no axis of reflection 
symmetry~\cite{Klaassen_2016}. For example, the five-fold rotationally symmetric 
edge-to-edge tiling by the convex hexagonal tile, where ``$A = D = 72^ \circ ,\;
B = E = 134^ \circ ,\;C = F = 154^ \circ ,\; a = b = c = d = e = f$" 
shown in Figure~\ref{fig05}(a), has $C_{5}$ symmetry. Additionally, the 
five-fold rotationally symmetric non-edge-to-edge tiling of the convex pentagonal 
tile generated by bisecting the convex hexagon, shown in Figure~\ref{fig05}(b), 
has $C_{5}$ symmetry. As described above, for $n \ge 4$ each $n$-fold, 
rotationally symmetric edge-to-edge tiling by a convex hexagonal tile that 
satisfies (\ref{eq2}) has $D_{n}$ symmetry. Conversely, the $n$-fold 
rotationally symmetric tilings by a convex pentagonal tile or a convex quadrangle 
generated from the convex hexagonal tile, which satisfies (\ref{eq2}) and is bisected 
by a dividing line that does not overlap with the axis of reflection 
symmetry, have $C_{n}$ symmetry. Therefore, the tilings of Figures~\ref{fig02}(b), 
\ref{fig02}(c), \ref{fig02}(d), and \ref{fig03}(b) have $C_{5}$ symmetry.

\begin{figure}[htb]
 \centering\includegraphics[width=14cm,clip]{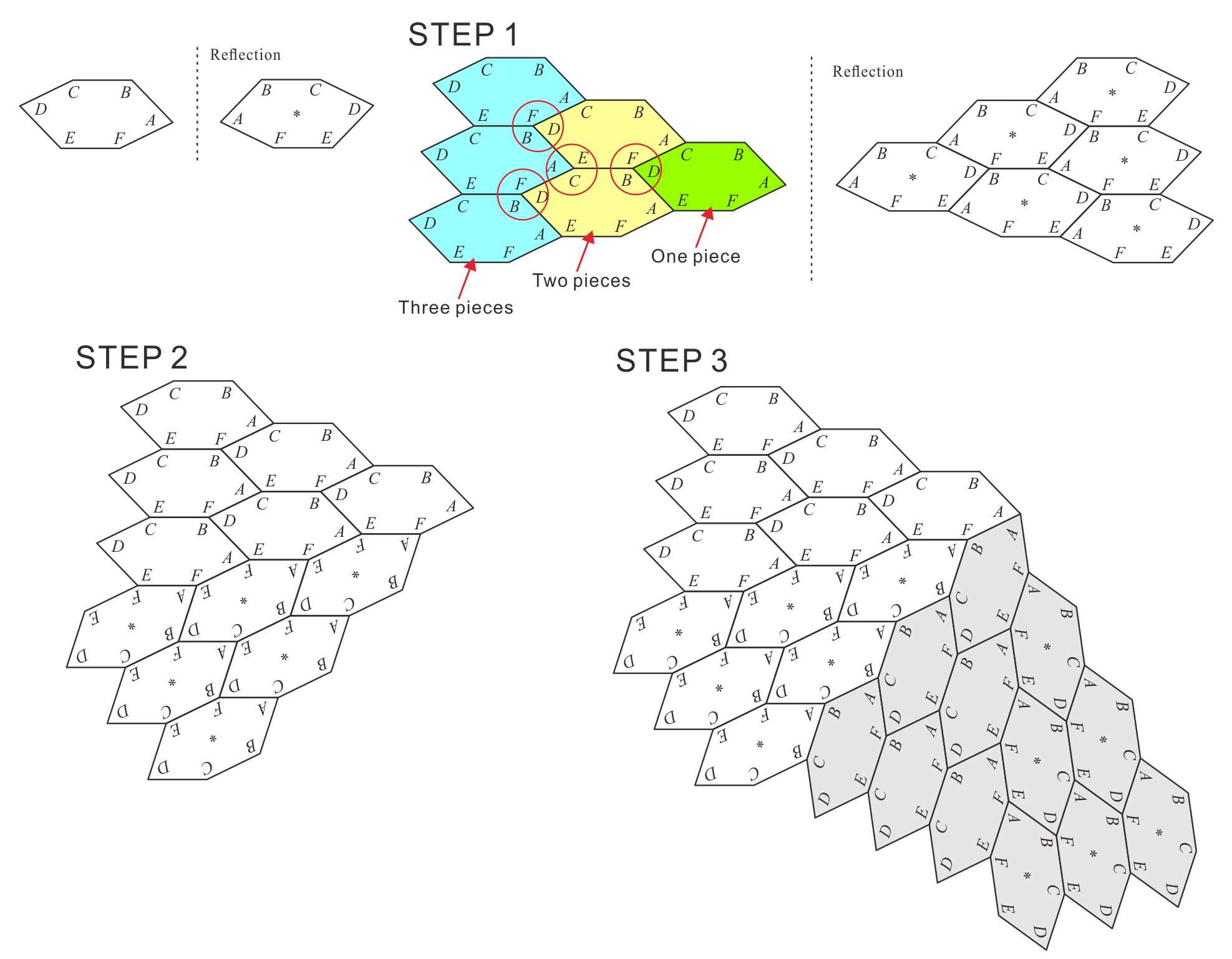} 
  \caption{{\small 
Formation method of rotationally symmetric edge-to-edge tiling with 
convex hexagonal tiles
} 
\label{fig06}
}
\end{figure}

Here, the formation of rotationally symmetric edge-to-edge tiling with 
convex hexagonal tiles is briefly explained. First, as shown in STEP 1 in 
Figure~\ref{fig06}, create a unit connecting the convex hexagonal tiles so that they 
form $n$-fold rotationally symmetric edge-to-edge tiling in one direction so 
that $B + D + F = 360^ \circ $ and $A + C + E = 360^ \circ $. The tiles can 
then be assembled in such a way as to increase the number of pieces from one 
to two to three, and so on, in order. Then create a similar unit with the 
reflected hexagons. Next, connect the two units created in STEP 1, so that 
$A + E + F = 360^ \circ $ as in STEP 2 in Figure~\ref{fig06}. 
Subsequently, take the unit from STEP 2 and rotate it by the value of the 
interior angle of vertex $A$. When the original unit and the rotated unit are 
arranged, so that $A + B + F = 360^ \circ $ and $A + B + C = 360^ \circ $ 
as shown in STEP 3 in Figure~\ref{fig06}, $\frac{2}{n}$ parts in the $n$-fold 
rotationally symmetric tiling can be formed. Then, by repeating this 
process as many times as necessary, an $n$-fold rotationally 
symmetric edge-to-edge tiling with convex hexagonal tiles can be formed. 
When the convex hexagons are bisected as depicted in Figure~\ref{fig02}, 
it will result in an $n$-fold rotationally symmetric non-edge-to-edge 
tiling with convex pentagonal tiles. Note that when the hexagons with 
the $D_{2}$ symmetry are used, as presented in Table~\ref{tab1}, a unit of 
reflected hexagons is not required. In this case, the units are arranged 
so that $A + B + E = 360^ \circ $ and $A + C + F = 360^ \circ $ as 
in STEP 2, and $A + B + F = 360^ \circ $, $A + C + F = 360^ \circ $, 
and $A + B + E = 360^ \circ $ as in STEP 3.

\renewcommand{\figurename}{{\small Figure.}}

%%%%%%%%%%%%%%%%%%%%%%%%%%%%%%%%%%%%%%%%%%%%%%%%%%%%%%%%%%%%%%%%%%%%%%
%%%%%%%%%%%%%%%%%%%%%%%%%%%%%%%%%%%%%%%%%%%%%%%%%%%%%%%%%%%%%%%%%%%%%%

\section{Rotationally symmetric tilings (tiling-like patterns) with a 
regular polygonal hole at the center}
\label{section3}

In \cite{Iliev_2018}, there are figures of rotationally symmetric tilings with a 
hole in the center of a regular hexagon, octagon, or 12-gon formed by using convex 
pentagons (or elements that can be regarded as convex hexagons). The 
regular hexagonal hole can be filled with convex pentagons. However, 
because the other holes cannot be filled with convex pentagons, they are 
not exactly tilings. The convex pentagons that can form rotationally 
symmetric tilings with a regular $m$-gonal hole at the center 
can be generated using the convex hexagons that satisfy (\ref{eq1}). 
The conditions of the convex hexagonal tiles that can be formed with a 
rotationally symmetric tiling with a regular $m$-gonal hole are expressed 
in (\ref{eq3}).

\begin{equation}
\label{eq3}
\left\{ {\begin{array}{l}
 A = D = \dfrac{720^ \circ }{m\strut}, \\ 
 B = C = E = F = 180^ \circ -\dfrac{A}{2} = 180^ \circ -\dfrac{360^ \circ}{m\strut}, \\
 a = b = c = d = e = f, \\ 
 \end{array}} \right.
\end{equation}

\noindent
where $m$ is an integer greater than or equal to five, because
 $0^\circ < A < 180^ \circ $. The convex hexagons that satisfy (\ref{eq3}) 
 are equilateral convex parallelohexagons with $D_{2}$ symmetry. Note that, 
 because ``$180^ \circ -\frac{360^ \circ}{m}$'' corresponds to one 
inner angle (interior angle) of a regular $m$-gon, the value of ``$A+F$" in (\ref{eq3}) 
is equal to the outer angle ($180^ \circ + \frac{360^ \circ}{m}$) 
of one vertex of a regular $m$-gon.

\begin{table}[htb]
 \begin{center}
{\small
\caption[Table 2]{
Interior angles of convex hexagons satisfying (\ref{eq3}) that 
can form the rotationally symmetric tilings with a regular $m$-gonal 
hole at the center
}
\label{tab2}
}
\begin{tabular}
{c| D{.}{.}{2} D{.}{.}{2} D{.}{.}{2} D{.}{.}{2} D{.}{.}{2} D{.}{.}{2} |c|c}
%{c|rrrrrr|c|c}
\hline
\raisebox{-1.50ex}[0cm][0cm]{$m$}& 
\multicolumn{6}{c|}{\shortstack{ Value of interior angle (degree) } } & 
\raisebox{-2.6ex}[0.7cm][0.5cm]{\shortstack{$n$ of \\Table~\ref{tab1}}}& 
\raisebox{-2.6ex}[0.7cm][0.5cm]{\shortstack{Figure \\number}} \\

 & 
\textit{A} & 
\textit{B}& 
\textit{C}& 
\textit{D}& 
\textit{E}& 
\textit{F}& 
  \\
\hline
5& 
144& 
108& 
108& 
144& 
108& 
108& 
& 
\ref{fig20} \\
\hline
6& 
120& 
120& 
120& 
120& 
120& 
120& 
3& 
\ref{fig21} \\
\hline
7& 
102.86 & 
128.57 & 
128.57 & 
102.86 & 
128.57 & 
128.57 & 
& 
\ref{fig22}\\
\hline
8& 
90& 
135& 
135& 
90& 
135& 
135& 
4& 
\ref{fig23}\\
\hline
9& 
80& 
140& 
140& 
80& 
140& 
140& 
& 
\ref{fig24}\\
\hline
10& 
72& 
144& 
144& 
72& 
144& 
144& 
5& 
\ref{fig25} \\
\hline
11& 
65.45 & 
147.27 & 
147.27 & 
65.45 & 
147.27 & 
147.27 & 
& 
 \\
\hline
12& 
60& 
150& 
150& 
60& 
150& 
150& 
6& 
\ref{fig26} \\
\hline
13& 
55.38 & 
152.31 & 
152.31 & 
55.38 & 
152.31 & 
152.31 & 
& 
 \\
\hline
14& 
51.43 & 
154.29 & 
154.29 & 
51.43 & 
154.29 & 
154.29 & 
7& 
\ref{fig27} \\
\hline
15& 
48& 
156& 
156& 
48& 
156& 
156& 
& 
 \\
\hline
16& 
45& 
157.5& 
157.5& 
45& 
157.5& 
157.5& 
8& 
\ref{fig28} \\
\hline
17& 
42.35 & 
158.82 & 
158.82 & 
42.35 & 
158.82 & 
158.82 & 
& 
 \\
\hline
18& 
40& 
160& 
160& 
40& 
160& 
160& 
9& 
 \\
\hline
19& 
37.89 & 
161.05 & 
161.05 & 
37.89 & 
161.05 & 
161.05 & 
& 
 \\
\hline
20& 
36& 
162& 
162& 
36& 
162& 
162& 
10& 
 \\
\hline
21& 
34.29 & 
162.86 & 
162.86 & 
34.29 & 
162.86 & 
162.86 & 
& 
 \\
\hline
22& 
32.73 & 
163.64 & 
163.64 & 
32.73 & 
163.64 & 
163.64 & 
11& 
 \\
\hline
23& 
31.30 & 
164.35 & 
164.35 & 
31.30 & 
164.35 & 
164.35 & 
& 
 \\
\hline
24& 
30& 
165& 
165& 
30& 
165& 
165& 
12& 
 \\
\hline
25& 
28.8& 
165.6& 
165.6& 
28.8& 
165.6& 
165.6& 
& 
 \\
\hline
...& 
...& 
...& 
...& 
...& 
...& 
...& 
& 
 \\
\hline
\end{tabular}

\end{center}

\end{table}

Table~\ref{tab2} presents some of the relationships between the interior angles of 
convex hexagons satisfying (\ref{eq3}) that can form the rotationally symmetric tilings 
with a regular $m$-gonal hole at the center. (For $m = 5\!-\!10, 12, 14, 16,$ 
tilings with a regular $m$-gonal hole at the center by convex hexagons and 
convex pentagons generated by bisecting them are drawn. For further details, 
Figures~\ref{fig20}--\ref{fig28}.) As mentioned above, if these elements are 
considered to be convex hexagons, the connection is edge-to-edge, and if they 
are considered to be convex pentagons created by bisecting the convex hexagon, 
the connection is non-edge-to-edge. These rotationally 
symmetric tilings with a regular $m$-gonal hole with $D_{m}$ symmetry 
at the center have $C_{m}$ symmetry. If convex hexagons satisfying 
(\ref{eq3}) have $m$ that is divisible by two, they are also convex hexagonal 
tiles that satisfy (\ref{eq2}). That is, the convex hexagonal tiles that satisfy 
(\ref{eq2}) can form a rotationally symmetric tiling with $C_{2n}$ symmetry with 
a regular $2n$-gonal hole with $D_{2n}$ symmetry at the center, and a rotationally 
symmetric tiling with $D_{n}$ symmetry. Then, the convex pentagonal tiles 
generated by bisecting convex hexagons that satisfy (\ref{eq2}) can form a 
rotationally symmetric tiling with $C_{2n}$ symmetry with a regular $2n$-gonal 
hole with $D_{2n}$ symmetry at the center, a rotationally symmetric tiling with 
$D_{n }$ symmetry, and rotationally symmetric tilings with $C_{n}$ 
symmetry.\footnote{ 
As shown in \cite{Klaassen_2016}, tilings with $D_{1}$ or 
$D_{2}$ symmetry are possible; however, the explanation for them has been 
omitted from this study.}

\renewcommand{\figurename}{{\small Figure.}}
\begin{figure}[htb]
 \centering\includegraphics[width=14cm,clip]{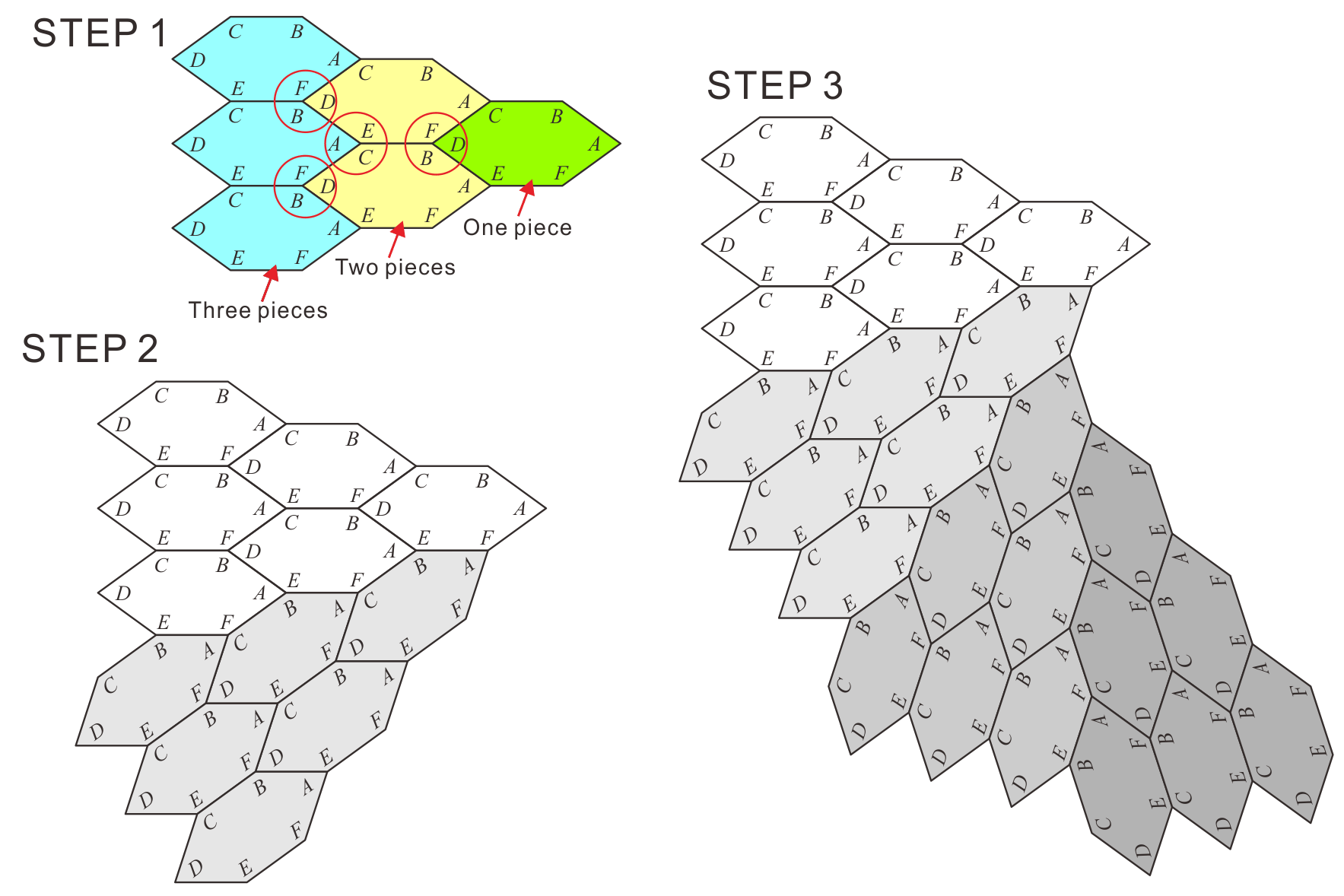} 
  \caption{{\small 
Formation method of rotationally symmetric tiling with a 
regular $m$-gonal hole at the center by convex hexagons 
in Table~\ref{tab2}
}
\label{fig07}
}
\end{figure}

Here, the formation of rotationally symmetric tiling with a regular $m$-gonal 
hole at the center by convex hexagons in Table 2 is briefly explained. 
First, as shown in STEP 1 in Figure~\ref{fig07}, create a unit connecting 
the convex hexagons satisfying (\ref{eq3}) in one direction so 
that $B + D + F = 360^ \circ $ and $A + C + E = 360^ \circ $. 
The hexagons can then be assembled in such a way as to increase the number 
of pieces from one to two to three, and so on, in order. Then, copy the unit 
in STEP 1, and rotate it by the half value of the interior  angle of vertex $A$. 
Subsequently, connect the two units of STEP 1, so that $A + B + E = 360^\circ$ 
and $A + C + F = 360^ \circ $ as in STEP 2 in Figure~\ref{fig07}. The series of two 
edges, \textit{AF}, of the unit in STEP 2 are edges of the contour of a regular 
$m$-gon. Then, by repeating this process as many times as necessary, a 
rotationally symmetric tiling with a regular $m$-gonal hole at the center with 
convex hexagons can be formed. STEP 3 in Figure~\ref{fig07} is $\frac{4}{10}$ 
parts of a regular 10-gon (The completed state is shown in Figure~\ref{fig25}). 
As shown in Table~\ref{tab2} and Figure~\ref{fig21}, the hexagon, $m = 6$, is a 
regular hexagon, so it is possible to fit the regular hexagon into the center hole.

%%%%%%%%%%%%%%%%%%%%%%%%%%%%%%%%%%%%%%%%%%%%%%%%%%%%%%%%%%%%%%%%%%%%%%
%%%%%%%%%%%%%%%%%%%%%%%%%%%%%%%%%%%%%%%%%%%%%%%%%%%%%%%%%%%%%%%%%%%%%%

\section{Rotationally symmetric tilings by an equilateral convex 
parallelohexagon with $C_{2}$ symmetry}
\label{section4}

The conditions of the equilateral convex parallelohexagons with $C_{2}$ 
symmetry are expressed in (\ref{eq4}).

\begin{equation}
\label{eq4}
\left\{ {\begin{array}{l}
 A = D \ne B, \\ 
 B = E \ne C, \\ 
 C = F \ne A, \\ 
 A + B + C = 360^ \circ , \\ 
 a = b = c = d = e = f. \\ 
 \end{array}} \right.
\end{equation}

\renewcommand{\figurename}{{\small Figure.}}
\begin{figure}[!h]
 \centering\includegraphics[width=15cm,clip]{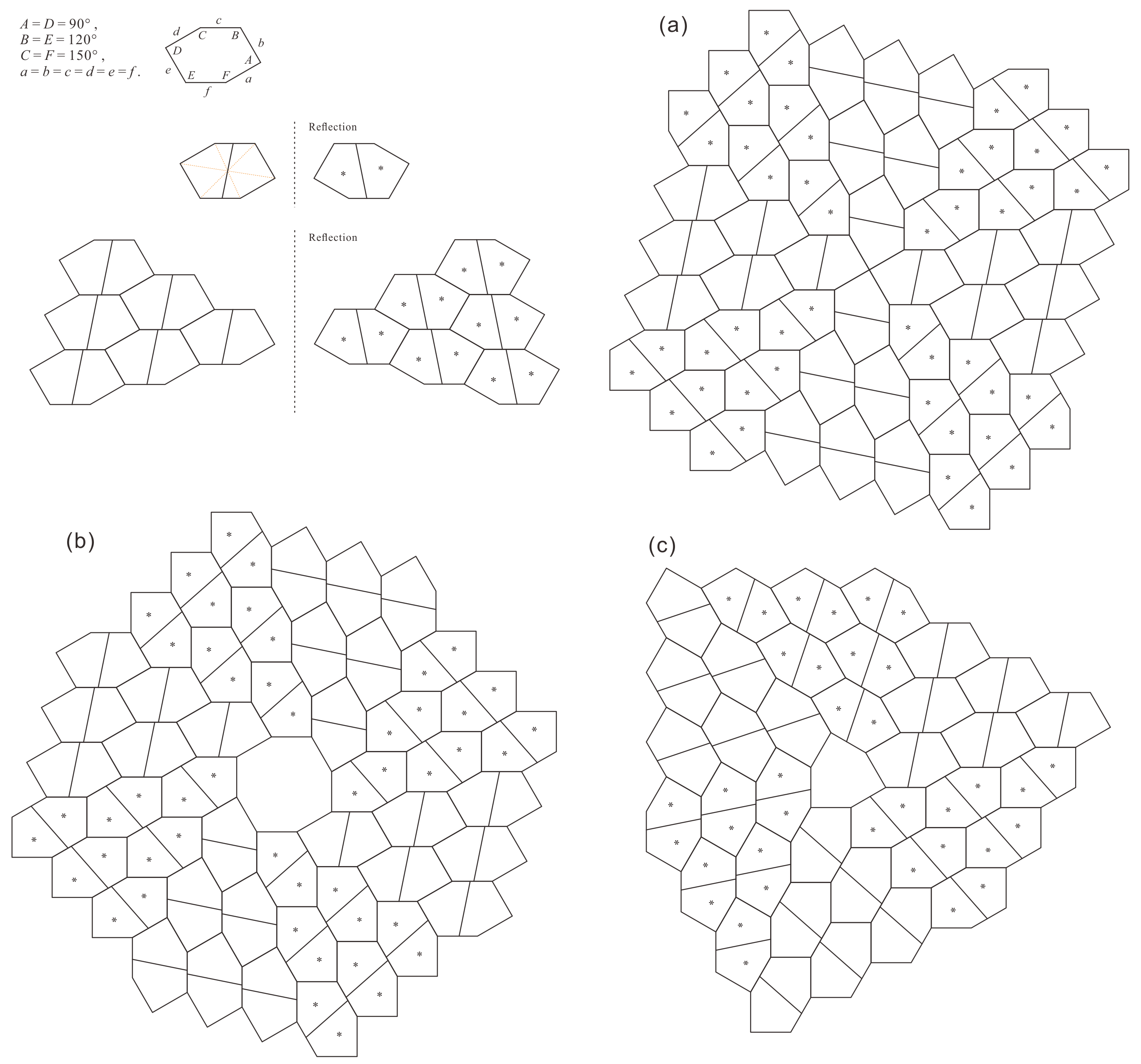} 
  \caption{{\small 
Rotationally symmetric tilings with $C_{3}$ or $C_{4}$ symmetry by 
a convex pentagon based on a convex hexagon with  $C_{2}$ symmetry 
that satisfies ``$A = D = 90^ \circ ,\;B = E = 120^ \circ ,\;C = F = 150^ \circ, \; 
a = b = c = d = e = f$")
}
\label{fig08}
}
\end{figure}

\renewcommand{\figurename}{{\small Figure.}}
\begin{figure}[!h]
 \centering\includegraphics[width=15cm,clip]{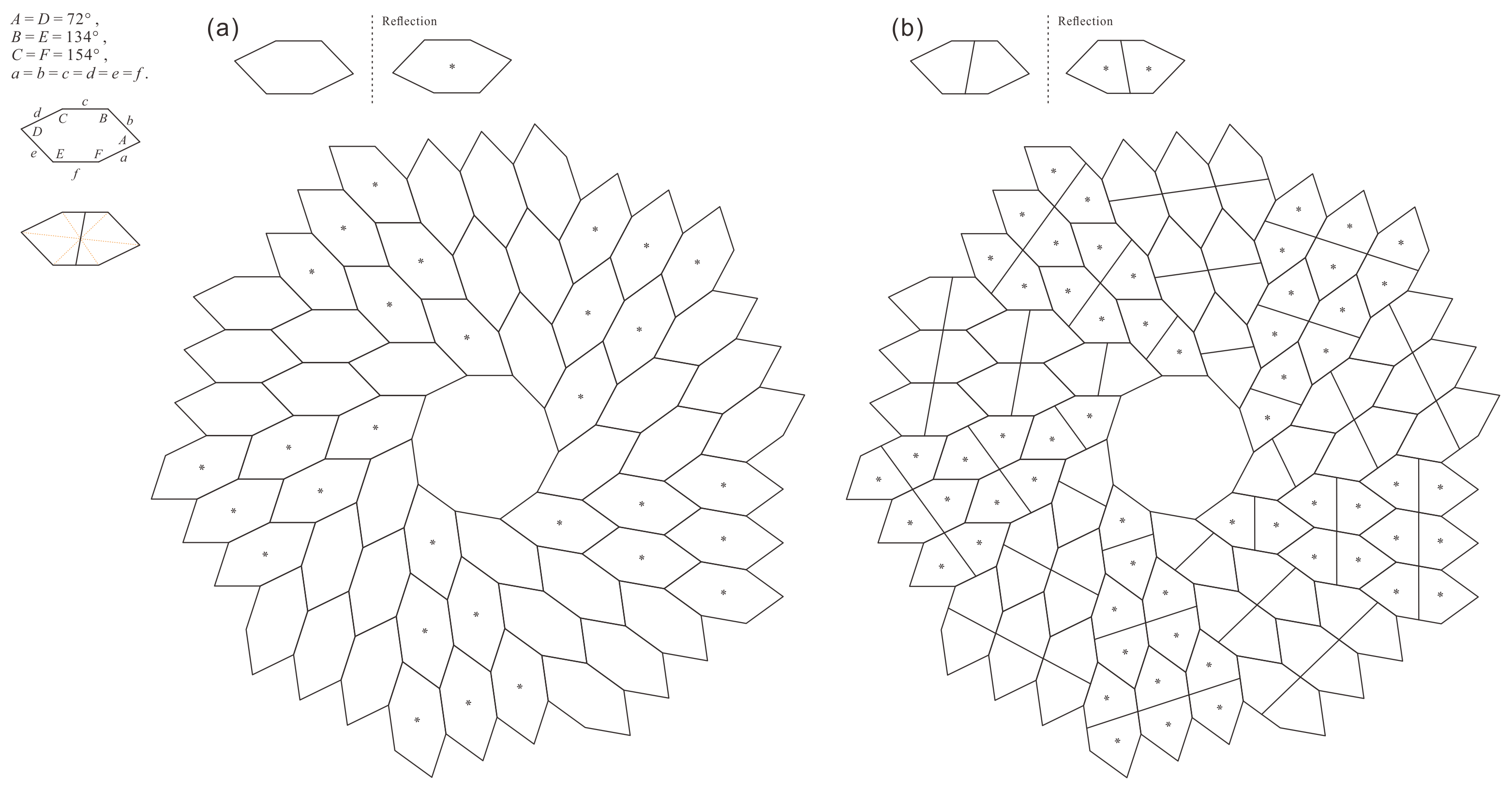} 
  \caption{{\small 
Rotationally symmetric tiling with $C_{5}$ symmetry with an equilateral 
convex 10-gonal hole with $D_{5}$ symmetry at the center using a convex 
hexagon with $C_{2}$ symmetry that satisfies ``$A = D = 72^ \circ ,\; 
B = E = 134^ \circ ,\;C = F = 154^ \circ ,\;a = b = c = d = e = f$," and the 
version of the tiling with convex pentagons based on the convex hexagons
}
\label{fig09}
}
\end{figure}

\renewcommand{\figurename}{{\small Figure.}}
\begin{figure}[!h]
 \centering\includegraphics[width=15cm,clip]{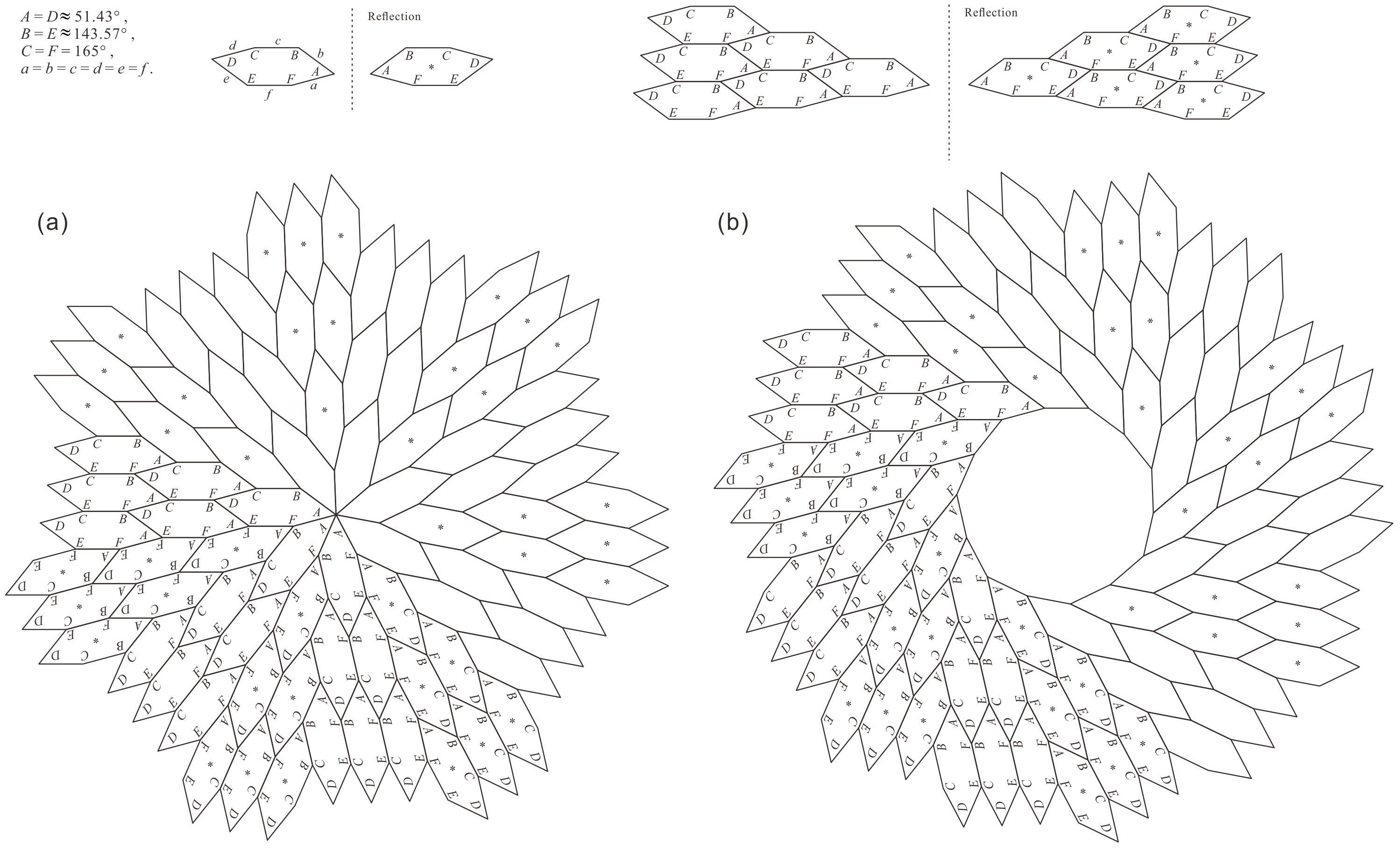} 
  \caption{{\small 
Rotationally symmetric tiling with $C_{7}$ symmetry and rotationally symmetric 
tiling with $C_{7}$ symmetry with an equilateral convex 14-gonal hole with 
$D_{7}$ symmetry at the center using an equilateral convex 
parallelohexagon with $C_{2}$ symmetry that satisfies 
``$A = D = \frac{360^ \circ }{7} \approx 51.43^ \circ , \; B = E \approx 143.57^ \circ, \;
C = F =165^ \circ, \;a = b = c = d = e = f$"
}
\label{fig10}
}
\end{figure}

\renewcommand{\figurename}{{\small Figure.}}
\begin{figure}[htb]
 \centering\includegraphics[width=15cm,clip]{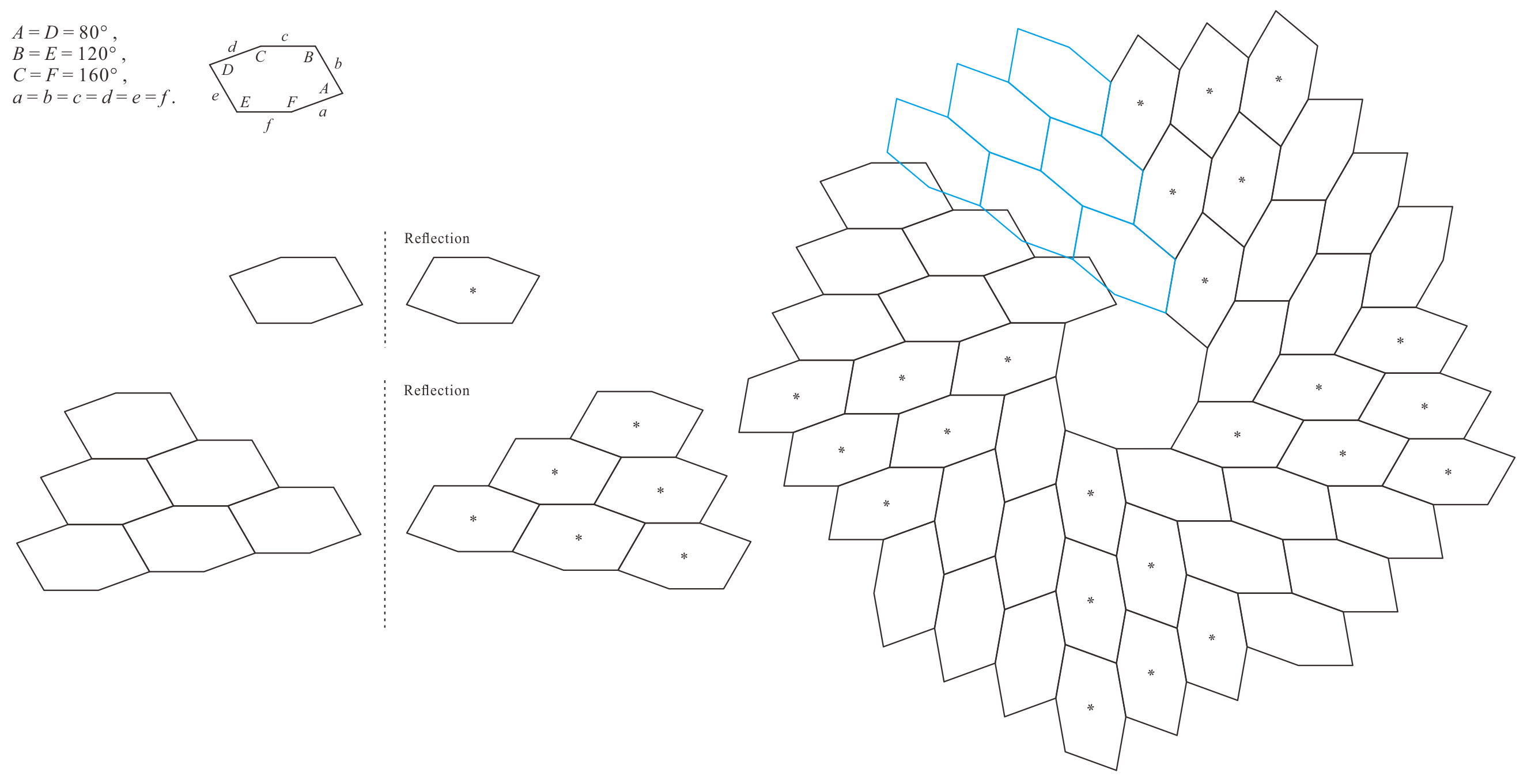} 
  \caption{{\small 
State in which a nonagonal hole cannot be formed using convex 
hexagons that satisfy "$A = D = \frac{720^ \circ }{9} = 80^ \circ ,\;B = E 
\ne C = F,\;A + B + C = 360^ \circ ,\;a = b = c = d = e = f$"
}
\label{fig11}
}
\end{figure}

As mentioned in Section \ref{section2}, the convex hexagonal tiles that satisfy 
(\ref{eq4}) and have interior angles of $\frac{360^ \circ }{n}$, where $n$ is an integer 
greater than or equal to three, can form rotationally symmetric tilings 
with $C_{n}$ symmetry.

In \cite{Iliev_2018}, Iliev presents some tilings using a convex pentagon that bisects 
an equilateral convex parallelohexagon with $C_{2}$ symmetry that satisfies 
``$A = D = 90^ \circ ,\;B = E = 120^ \circ ,\;C = F = 150^ \circ ,\; 
a = b = c = d = e = f$." They are a rotationally symmetric tiling with $C_{4}$ 
symmetry, as shown in Figure~\ref{fig08}(a), a rotationally symmetric tiling with 
$C_{4}$ symmetry with an equilateral convex octagonal hole with $D_{4}$ symmetry 
at the center, as shown in Figure~\ref{fig08}(b), and a rotationally symmetric tiling 
with $C_{3}$ symmetry with an equilateral convex hexagonal hole with $D_{3}$ 
symmetry at the center, as shown in Figure~\ref{fig08}(c). Notably, those tilings can 
be formed with convex octagonal and hexagonal holes that are equilateral but 
not regular polygons at the center, as shown in Figures~\ref{fig08}(b) and \ref{fig08}(c). 
Thus, the convex hexagonal tiles that satisfy (\ref{eq4}) and have interior angles of 
$\frac{360^ \circ }{n}$ can form rotationally symmetric tilings with $C_{n}$ symmetry 
with an equilateral convex 2$n$-gonal hole with $D_{n}$ symmetry at the center. 
For example, Figure~\ref{fig09}(a) presents a rotationally symmetric tiling with 
$C_{5}$ symmetry with an equilateral convex 10-gonal hole with $D_{5}$ symmetry 
at the center using an equilateral convex parallelohexagon with $C_{2}$ symmetry 
that satisfies ``$A = D = \frac{360^ \circ }{5} = 72^ \circ ,\; B = E = 134^ \circ ,\;
C = F = 154^ \circ ,\;a = b = c = d = e = f$," and Figure~\ref{fig09}(b) is a 
version of the tiling with convex pentagons bisecting the convex hexagons. 
Note that the convex hexagon in Figure~\ref{fig09} is the same as the convex 
hexagon in Figure~\ref{fig05}. Then, Figure~\ref{fig10}(a) shows a rotationally 
symmetric tiling with $C_{7}$ symmetry, and Figure~\ref{fig10}(b) shows a 
rotationally symmetric tiling with $C_{7}$ symmetry with an equilateral convex 
14-gonal hole with $D_{7}$ symmetry at the center, using an equilateral convex 
parallelohexagon with $C_{2}$ symmetry that satisfies 
``$A = D = \frac{360^ \circ }{7} \approx 51.43^ \circ, \; B = E \approx 143.57^ \circ, \;
C = F =165^ \circ, \; a = b = c = d = e = f$."

\renewcommand{\figurename}{{\small Figure.}}
\begin{figure}[!h]
 \centering\includegraphics[width=15cm,clip]{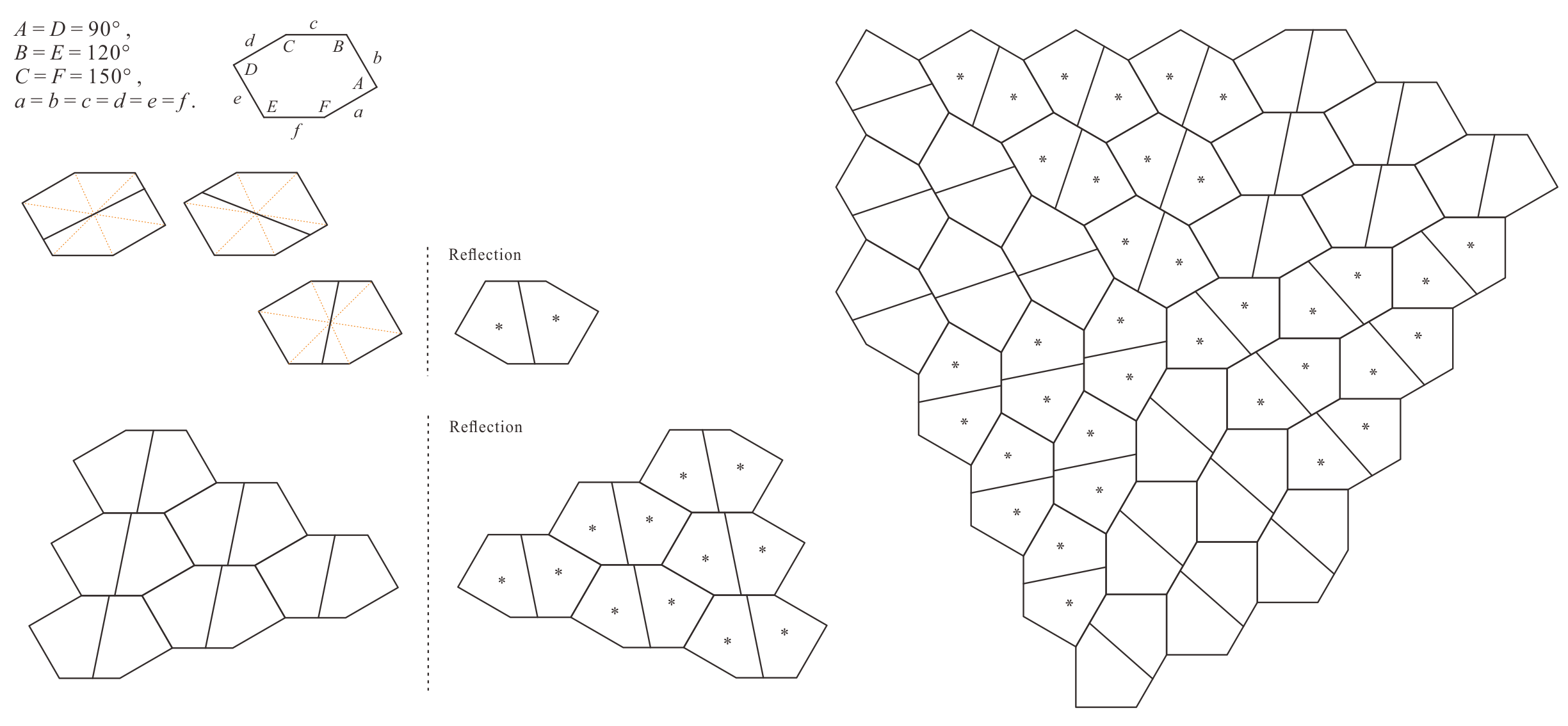} 
  \caption{{\small 
Three-fold rotationally symmetric tiling by a convex pentagon based 
on a convex hexagon that satisfies ``$A = D = 90^ \circ ,\;B = E = 120^ \circ ,\;
C = F = 150^ \circ, \; a = b = c = d = e = f$"
}
\label{fig12}
}
\end{figure}

Let us supplement the reason why the number of edges of the central polygonal 
hole, i.e., the equilateral convex polygonal hole at the center that is formed by 
an equilateral convex parallelohexagon with $C_{2}$ symmetry as described 
above, is even. For example, when $m = 5$ in Figure~\ref{fig20} (the five-fold 
rotationally symmetric tilings with a regular pentagonal hole by a convex hexagonal 
tile that satisfies (\ref{eq3}) and $A = D = \frac{720^ \circ }{5} = 144^ \circ )$, it can 
be observed that the five units created in STEP 1 in Figure~\ref{fig07} are used. 
Conversely, when forming a rotationally symmetric tiling by an equilateral convex 
parallelohexagon with the $C_{2}$ symmetry, as shown in Figures~\ref{fig05}, 
\ref{fig06}, \ref{fig08}, \ref{fig09} and \ref{fig10}, the units with hexagons and 
units with reflected hexagons should be connected alternately. Therefore, if 
the number of edges of the central polygonal hole is odd, the convex hexagons 
with $C_{2}$ symmetry cannot form the central polygon shown in Figure~\ref{fig11} 
(when trying to form a nonagonal hole by an equilateral convex parallelohexagon 
with $C_{2}$ symmetry that satisfies $A = D = \frac{720^ \circ }{9} = 80^ \circ )$. 
Because an equilateral convex $2n$-gon with $D_{n }$ symmetry has two types of 
inner angles (interior angles) in vertices, there are two types of outer angles in vertices. 
If $A = D = \frac{360^ \circ }{n}$ in (\ref{eq4}), values of ``$A + B,\;A + F$" are equal 
to the two types of outer angles (see Figure~\ref{fig10}(b)).

\begin{figure}[!h]
 \centering\includegraphics[width=15cm,clip]{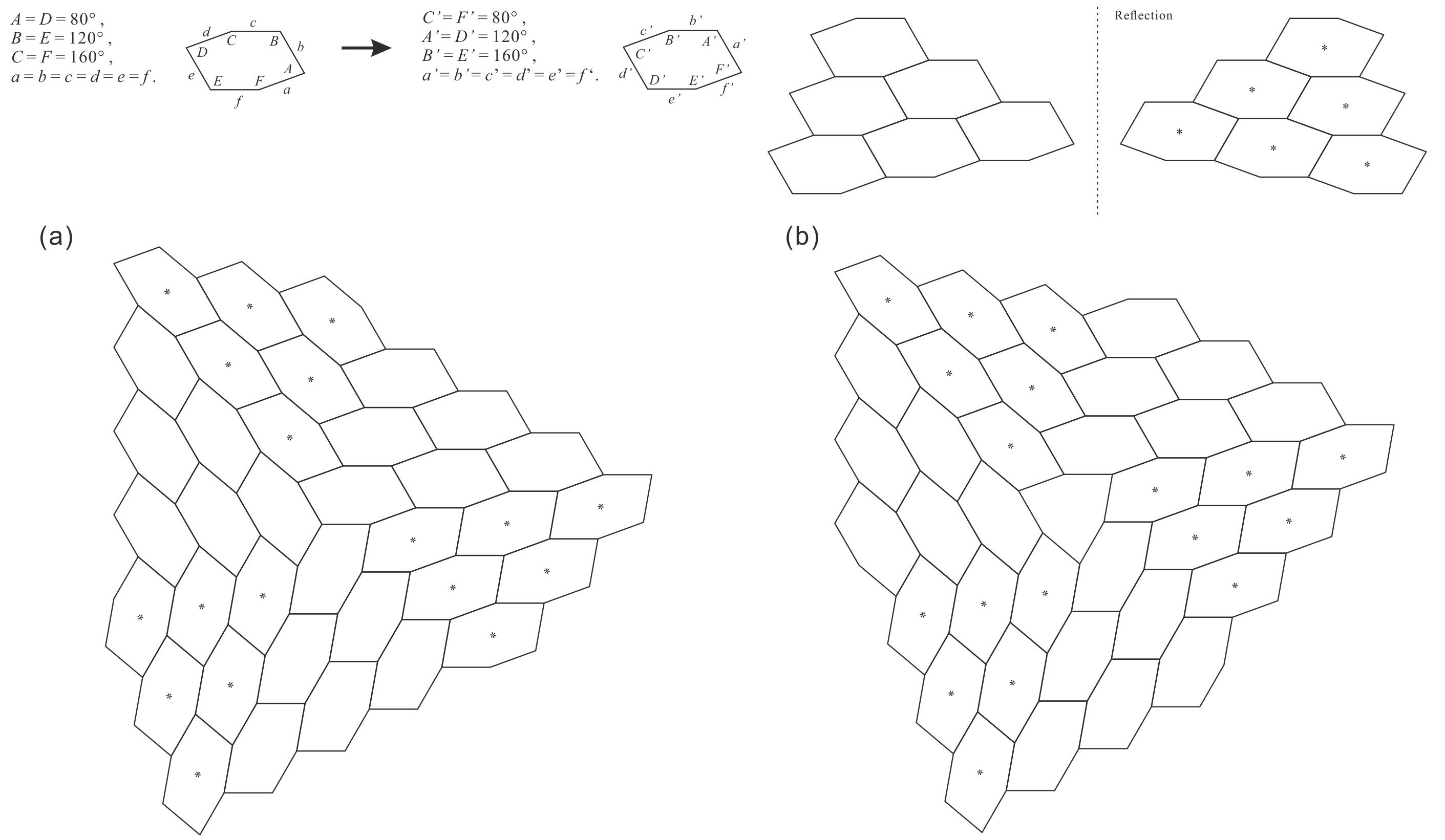} 
  \caption{{\small 
Three-fold rotationally symmetric tiling, and rotationally symmetric 
tiling with an equilateral convex hexagonal hole at the center by a convex 
hexagon that satisfies ``$A = D = 80^ \circ ,\;B = E = 120^ \circ ,\;
C = F = 160^ \circ ,\;a = b = c = d = e = f$"
}
\label{fig13}
}
\end{figure}

Figure~\ref{fig12} shows a rotationally symmetric tiling with $C_{3}$ symmetry formed by  
a convex hexagonal tile that satisfies ``$A = D = 90^ \circ ,\;B = E = 120^ \circ ,\;
C = F = 150^ \circ ,\;a = b = c = d = e = f$" in Figure 8. Then, Figure~\ref{fig13}(a) 
shows a three-fold rotationally symmetric tiling formed by the convex hexagonal 
tile that satisfies ``$A = D = 80^ \circ ,\;B = E = 120^ \circ ,\;C = F = 160^ \circ ,\;
a = b = c = d = e = f$" shown in Figure~\ref{fig11}. As shown in Figure~\ref{fig13}(b), 
the convex hexagons can also form a rotationally symmetric tiling with $C_{3}$ 
symmetry with an equilateral convex hexagonal hole at the center. This arrangement 
indicates that the convex hexagons that satisfy (\ref{eq4}), where one of the vertices 
$A$, $B$, or $C$ of which is equal to $\frac{360^ \circ }{n}$, can form an $n$-fold 
rotationally symmetric tiling, and a tiling with $C_{n}$ symmetry with an equilateral 
convex $2n$-gonal hole with $D_{n}$ symmetry at the center. (It is a matter of 
starting the vertex symbol from somewhere. As shown in Figure~\ref{fig13}, if 
``$A \to F',\;B \to A',\;\ldots $" is replaced, then $A' = D' = \frac{360^ \circ }{3} = 120^ \circ $, 
which corresponds to (\ref{eq1}).)

\begin{figure}[htbp]
 \centering\includegraphics[width=15cm,clip]{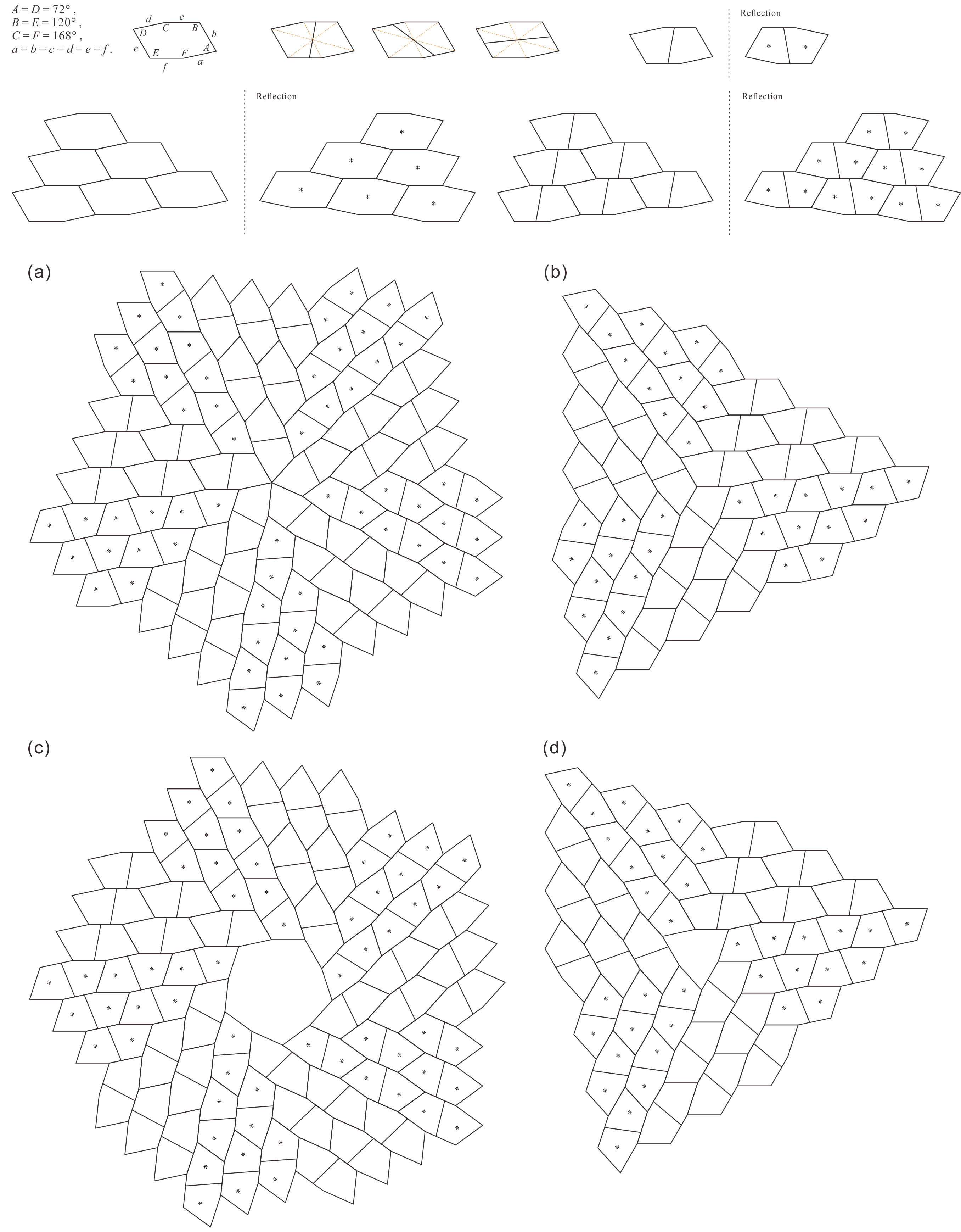} 
  \caption{{\small 
Five-fold or three-fold rotationally symmetric tilings, and 
rotationally symmetric tilings with an equilateral convex 10-gonal or 
hexagonal hole at the center by a convex pentagon based on a convex 
hexagon that satisfies ``$A = D = 72^ \circ ,\; B = E = 120^ \circ ,\; 
C = F = 168^ \circ ,\; a = b = c = d = e = f$" }
\label{fig14}
}
\end{figure}

\begin{figure}[!h]
 \centering\includegraphics[width=14cm,clip]{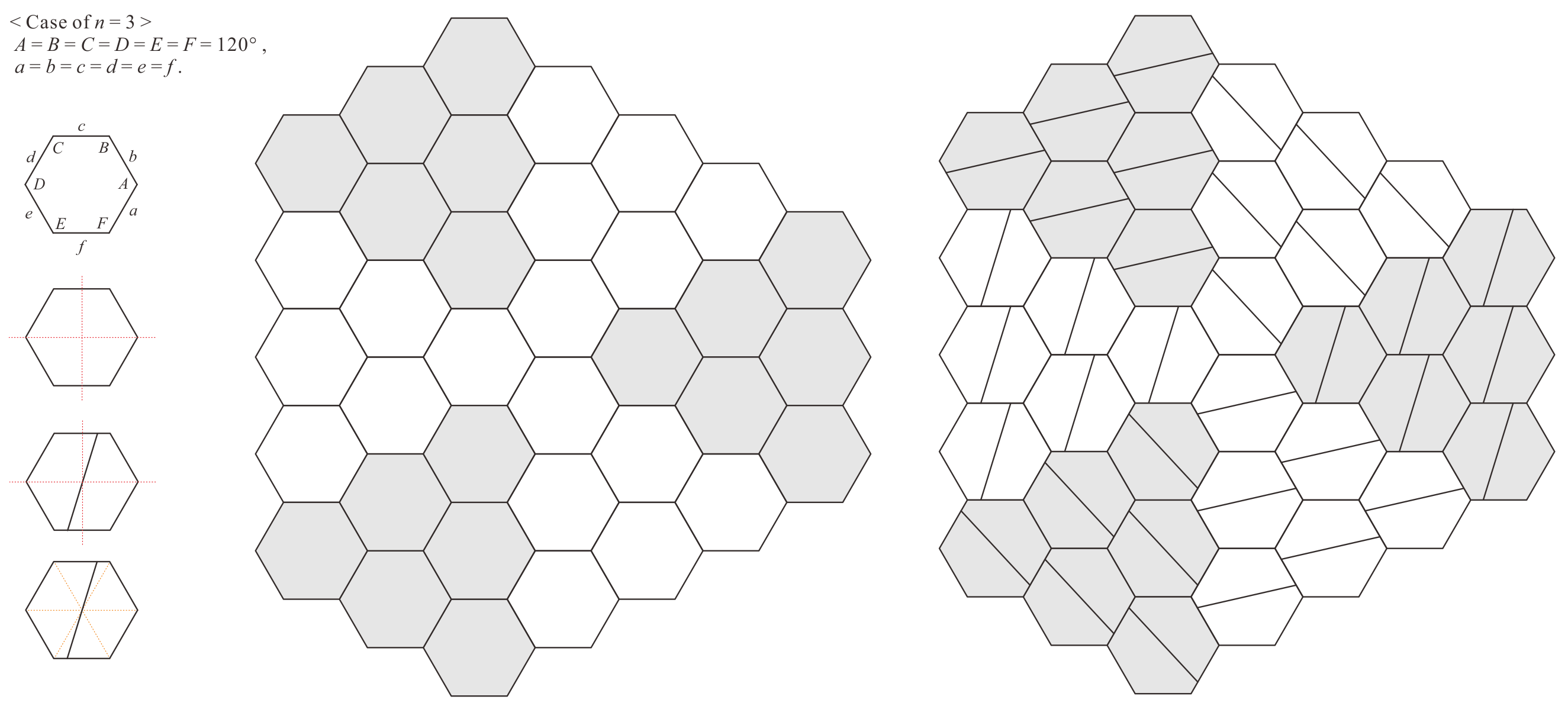} 
  \caption{{\small 
Three-fold rotationally symmetric tilings by a convex hexagon 
and a convex pentagon
}
\label{fig15}
}
\end{figure}

\begin{figure}[!h]
 \centering\includegraphics[width=14cm,clip]{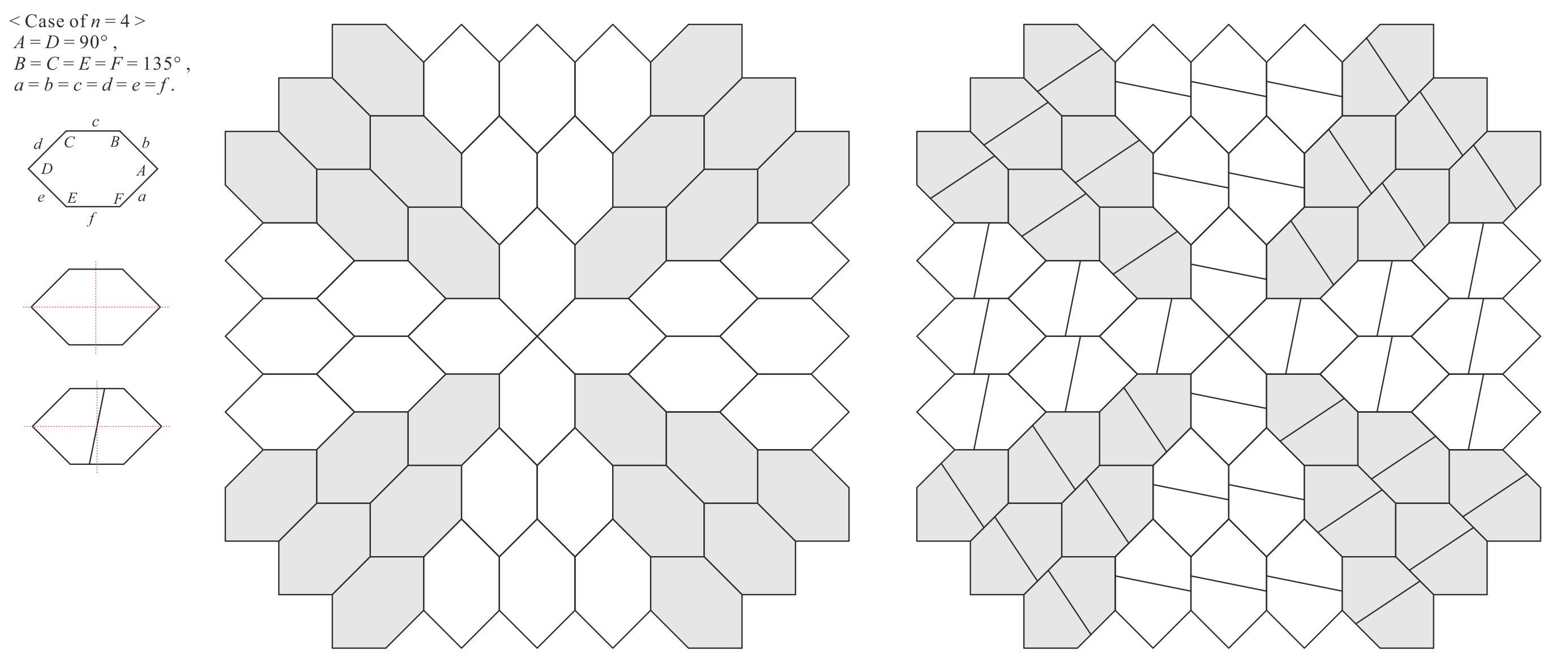} 
  \caption{{\small 
Four-fold rotationally symmetric tilings by a convex hexagon 
and a convex pentagon
}
\label{fig16}
}
\end{figure}

\begin{figure}[!h]
 \centering\includegraphics[width=14cm,clip]{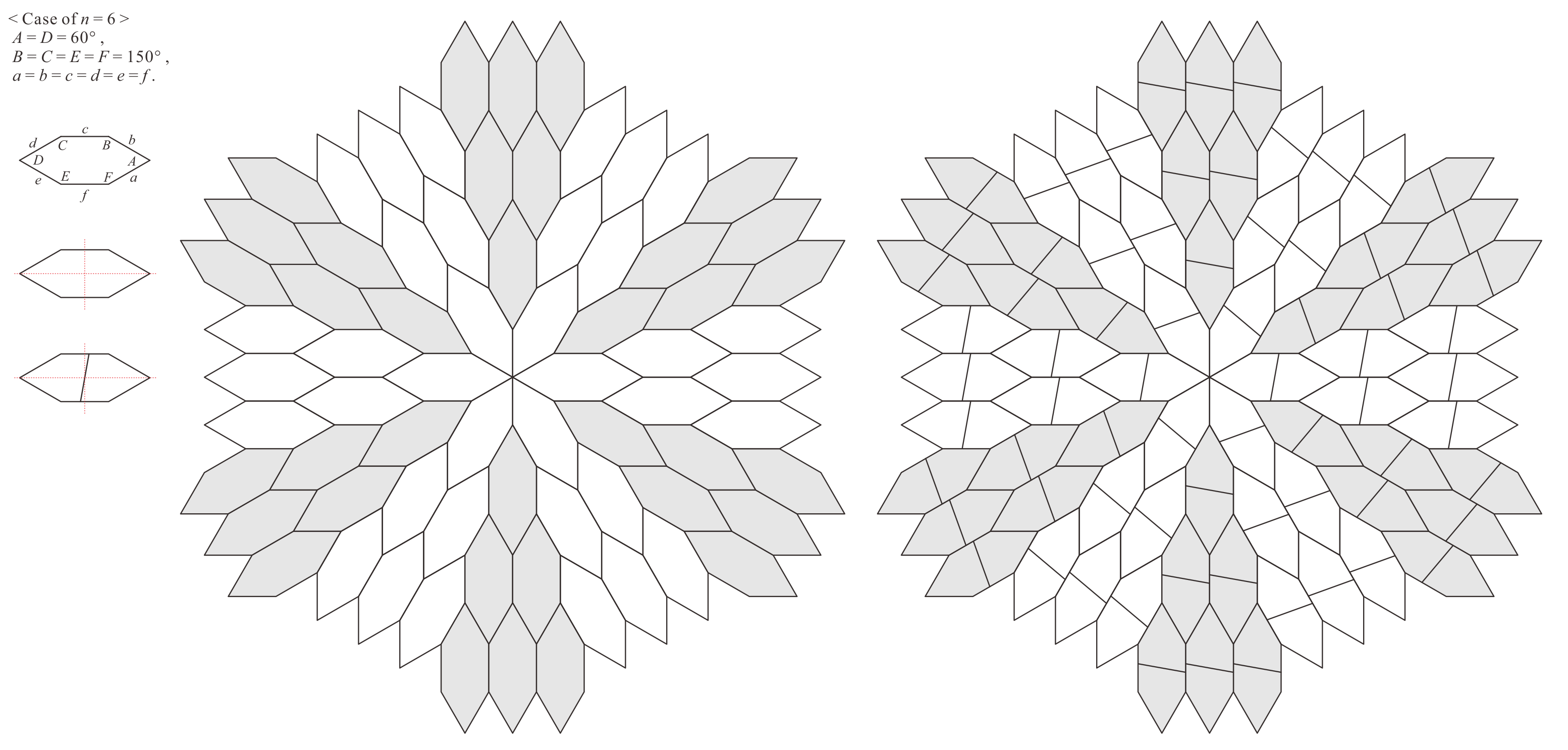} 
  \caption{{\small 
Six-fold rotationally symmetric tilings by a convex hexagon 
and a convex pentagon
}
\label{fig17}
}
\end{figure}

\begin{figure}[!h]
 \centering\includegraphics[width=14cm,clip]{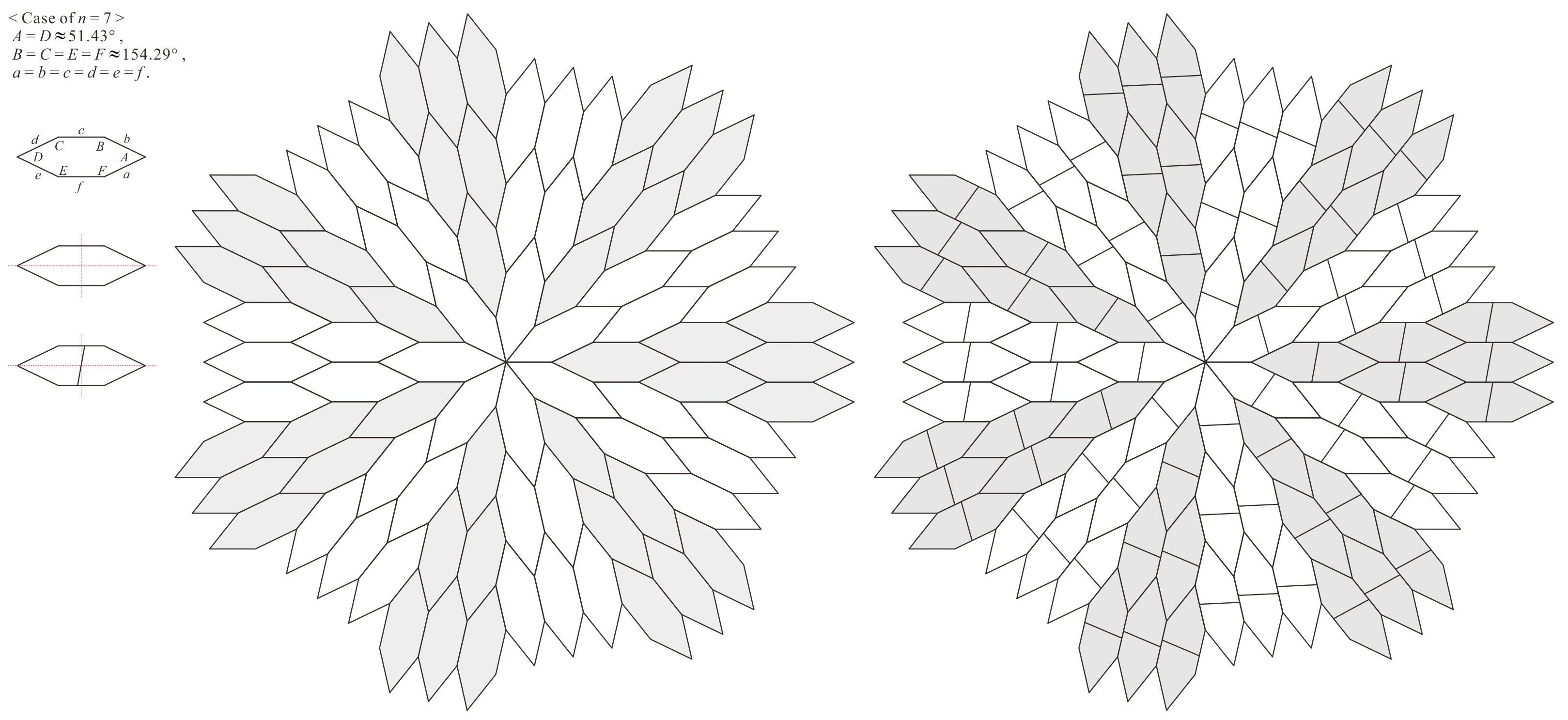} 
  \caption{{\small 
Seven-fold rotationally symmetric tilings by a convex hexagon 
and a convex pentagon
}
\label{fig18}
}
\end{figure}

\begin{figure}[htbp]
 \centering\includegraphics[width=15cm,clip]{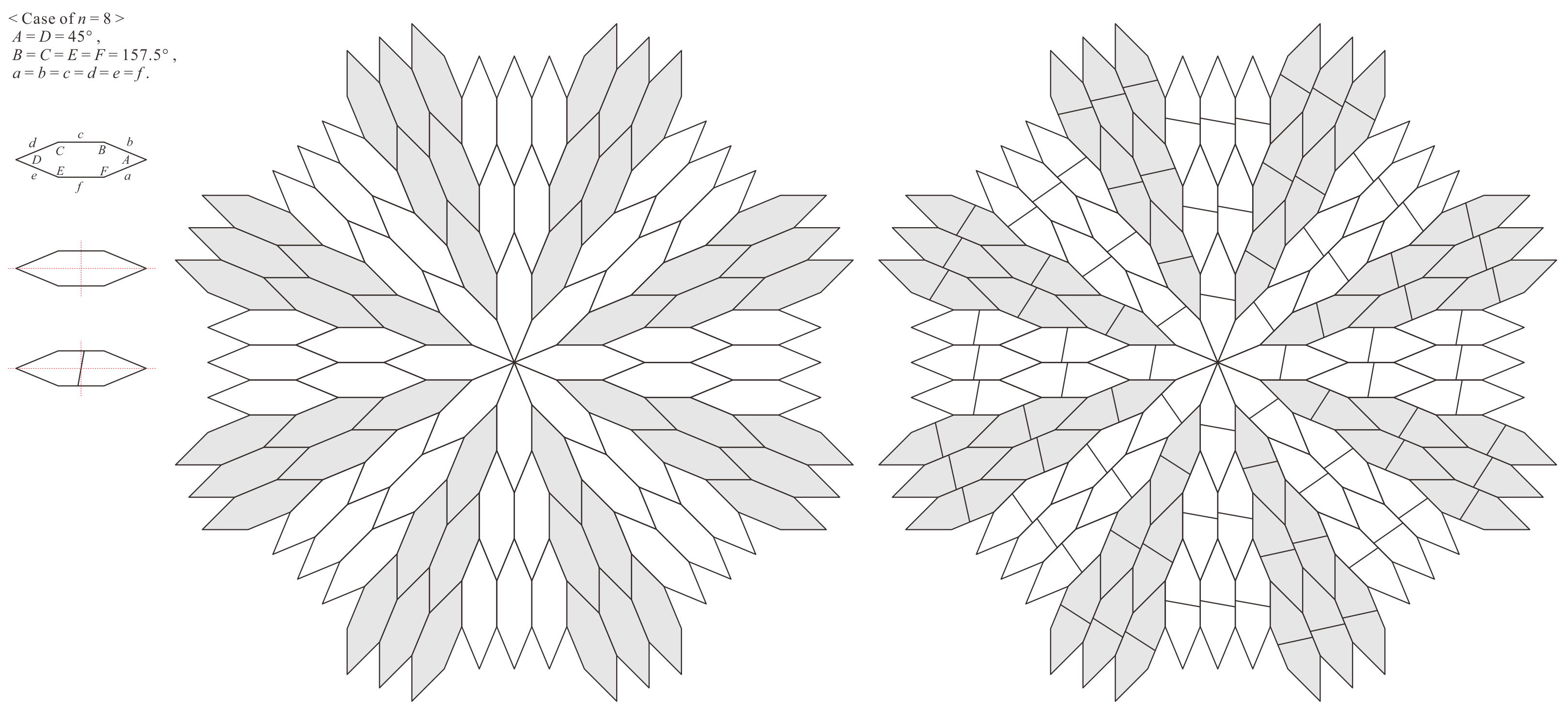} 
  \caption{{\small 
Eight-fold rotationally symmetric tilings by a convex hexagon 
and a convex pentagon
}
\label{fig19}
}
\end{figure}

\begin{figure}[htbp]
 \centering\includegraphics[width=11.5cm,clip]{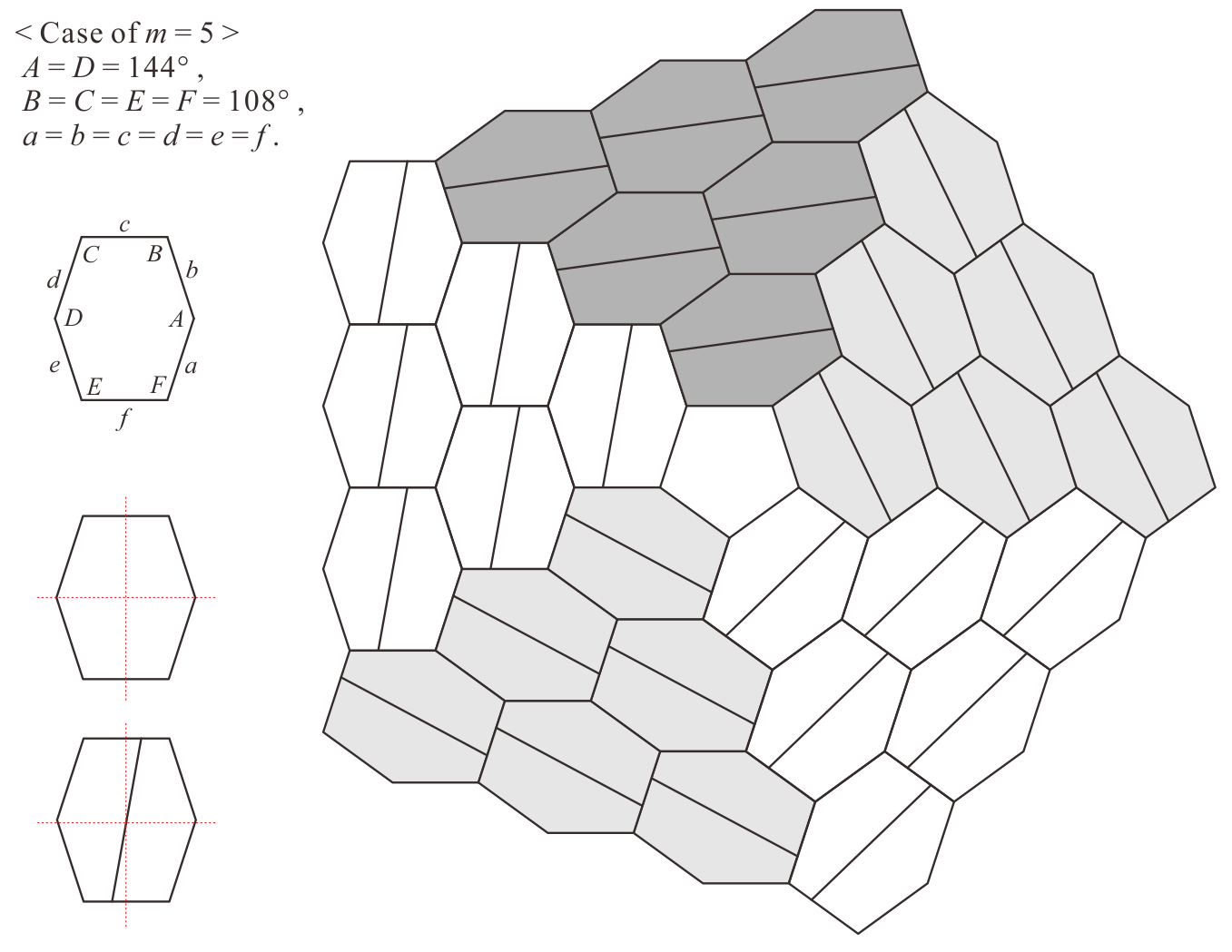} 
  \caption{{\small 
Rotationally symmetric tiling with $C_{5}$ symmetry with a regular 
convex pentagonal hole at the center by a convex pentagon
}
\label{fig20}
}
\end{figure}

\begin{figure}[htbp]
 \centering\includegraphics[width=11.5cm,clip]{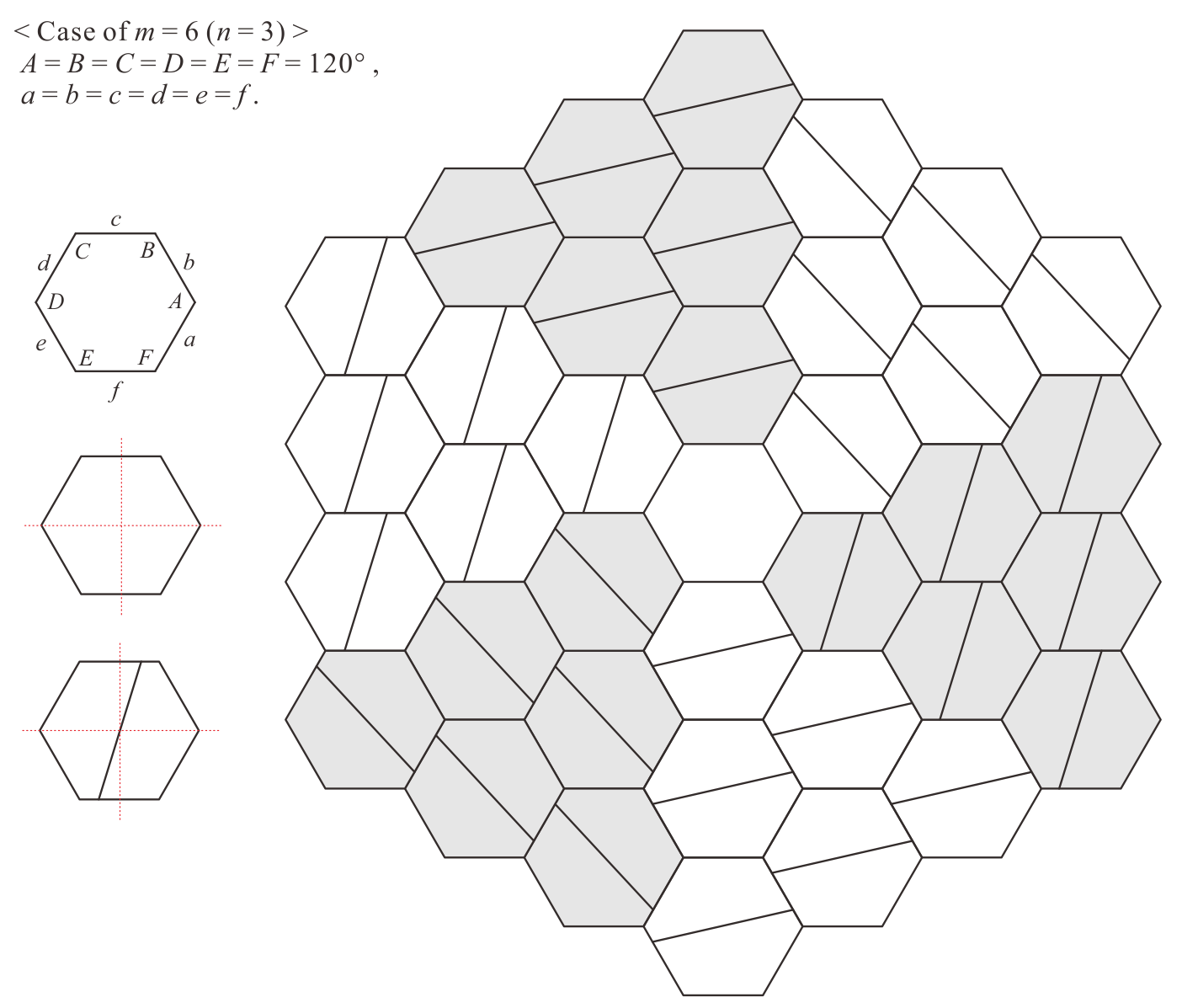} 
  \caption{{\small 
Rotationally symmetric tiling with $C_{6}$ symmetry with a regular 
convex hexagonal hole at the center by a convex pentagon
}
\label{fig21}
}
\end{figure}

\begin{figure}[htbp]
 \centering\includegraphics[width=11.5cm,clip]{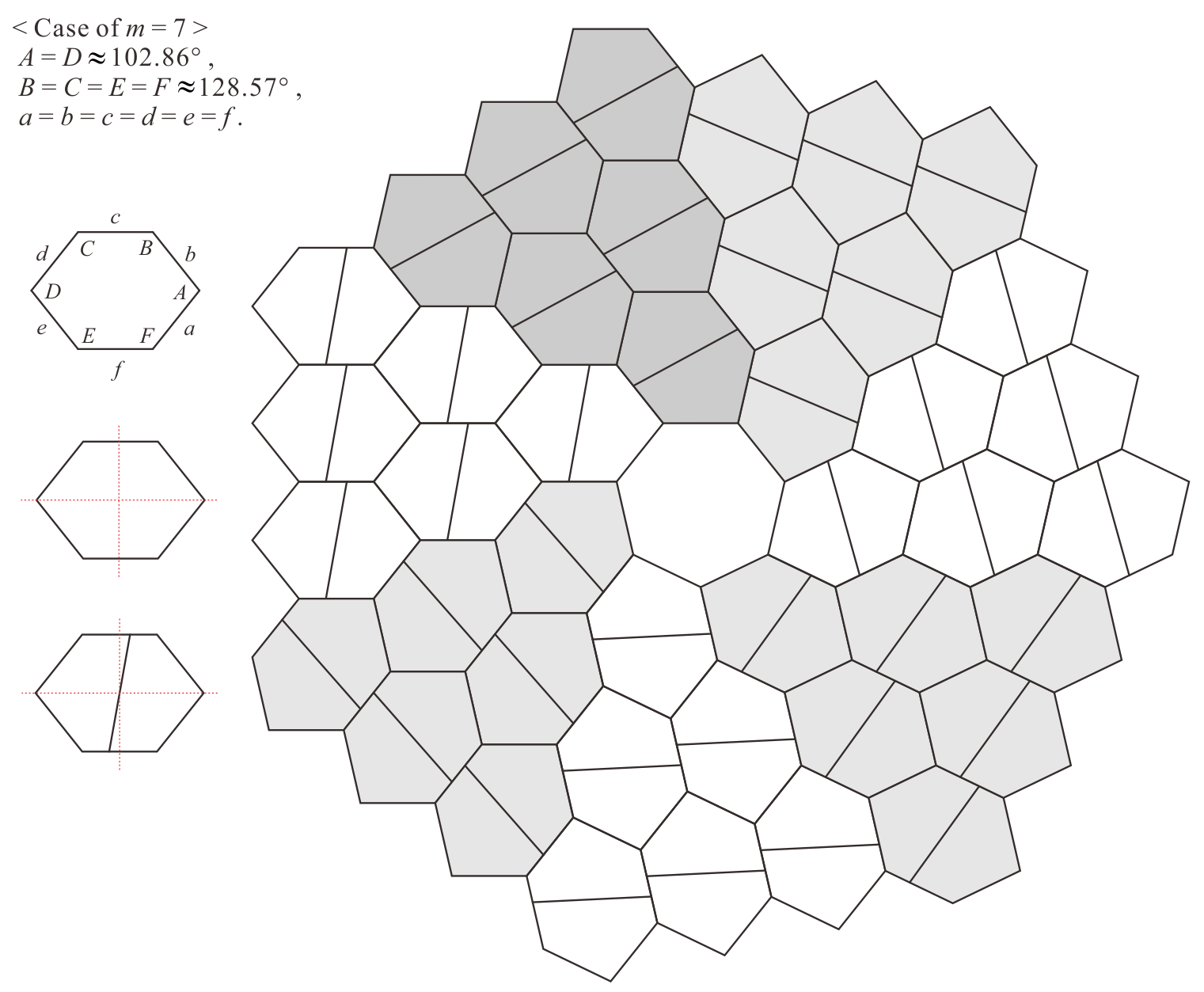} 
  \caption{{\small 
Rotationally symmetric tiling with $C_{7}$ symmetry with a regular 
convex heptagonal hole at the center by a convex pentagon
}
\label{fig22}
}
\end{figure}

\begin{figure}[htbp]
 \centering\includegraphics[width=11.5cm,clip]{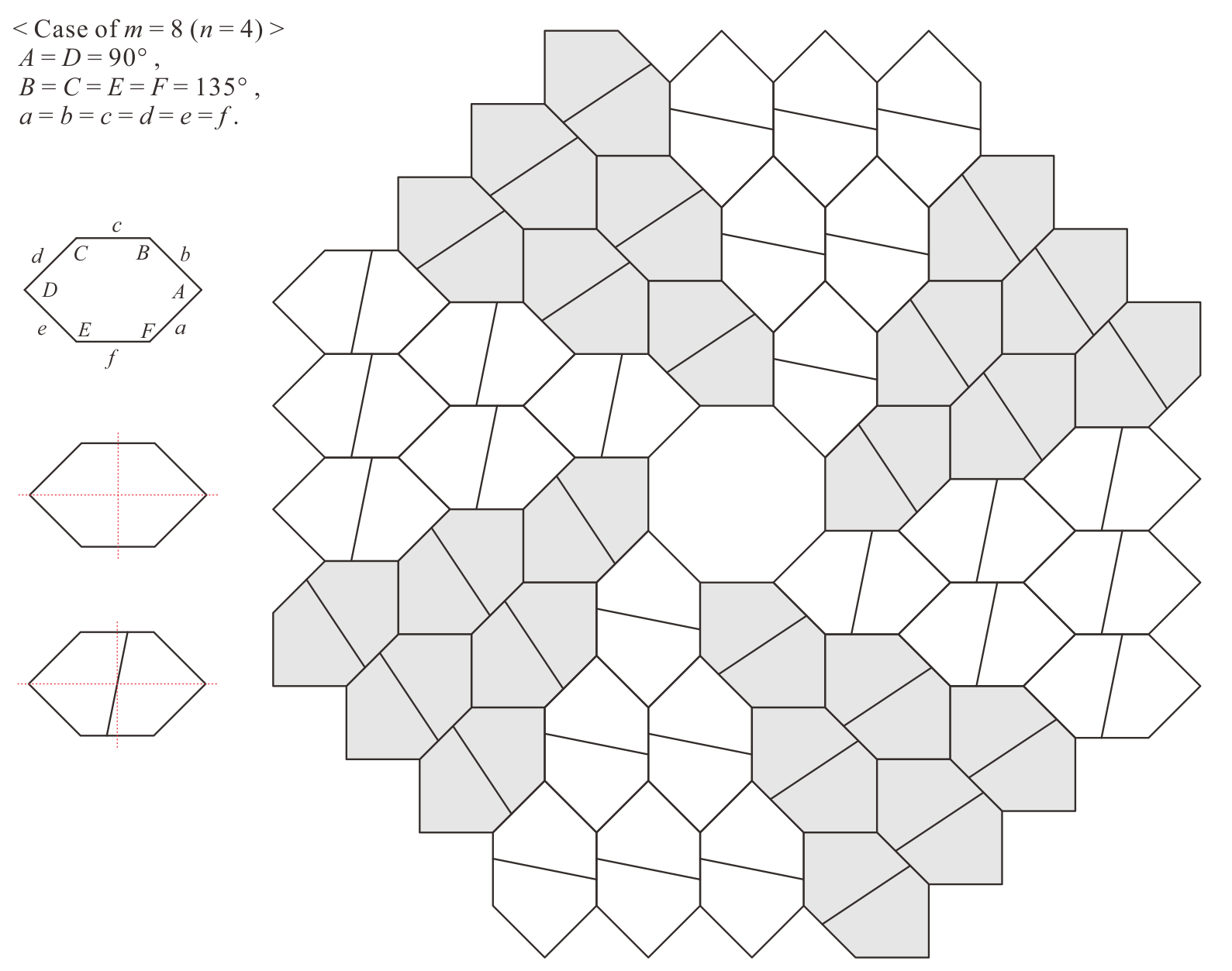} 
  \caption{{\small 
Rotationally symmetric tiling with $C_{8}$ symmetry with a regular 
convex octagonal hole at the center by a convex pentagon
}
\label{fig23}
}
\end{figure}

\begin{figure}[htbp]
 \centering\includegraphics[width=11.5cm,clip]{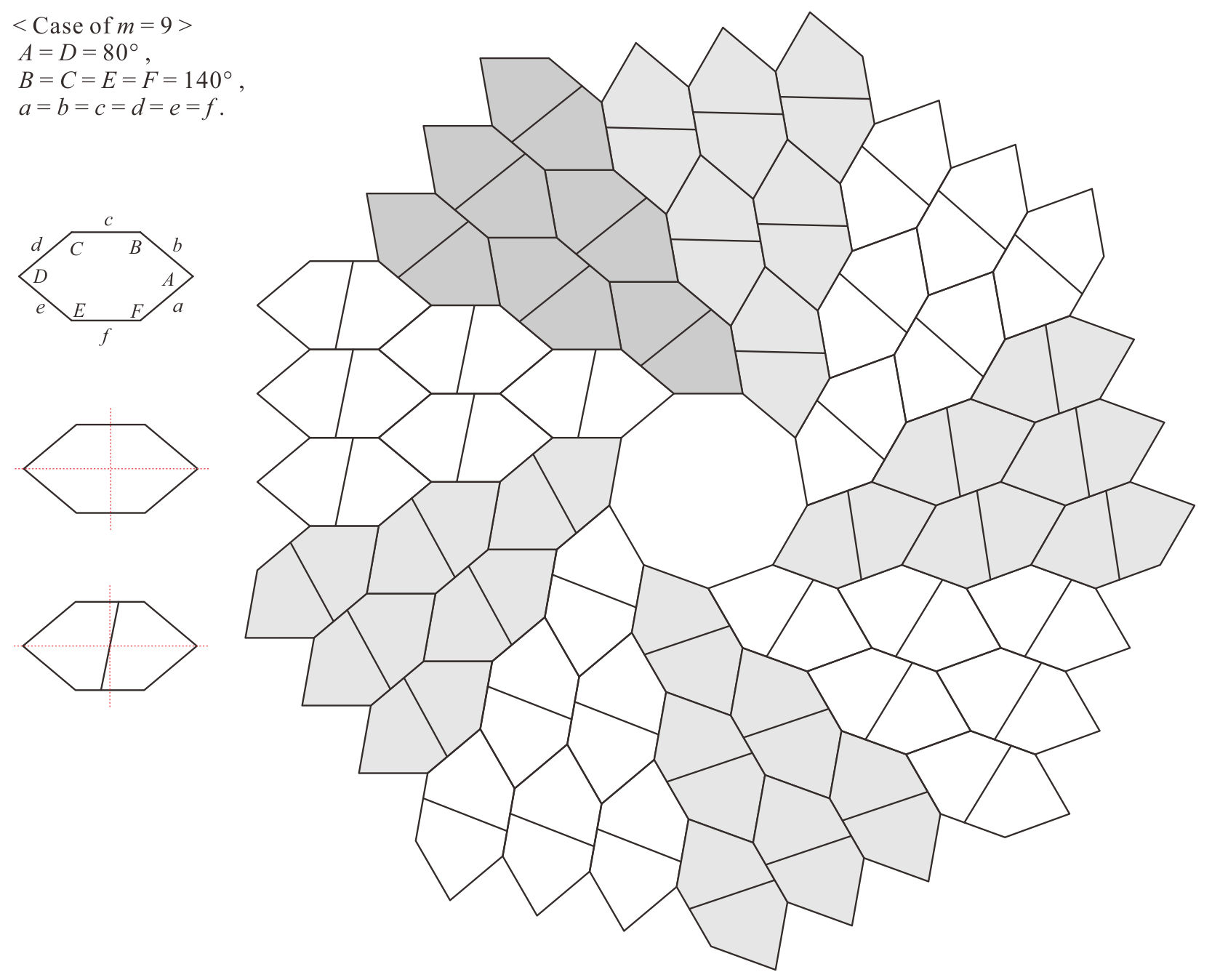} 
  \caption{{\small 
Rotationally symmetric tiling with $C_{9}$ symmetry with a regular 
convex nonagonal hole at the center by a convex pentagon
}
\label{fig24}
}
\end{figure}

\begin{figure}[htbp]
 \centering\includegraphics[width=11.5cm,clip]{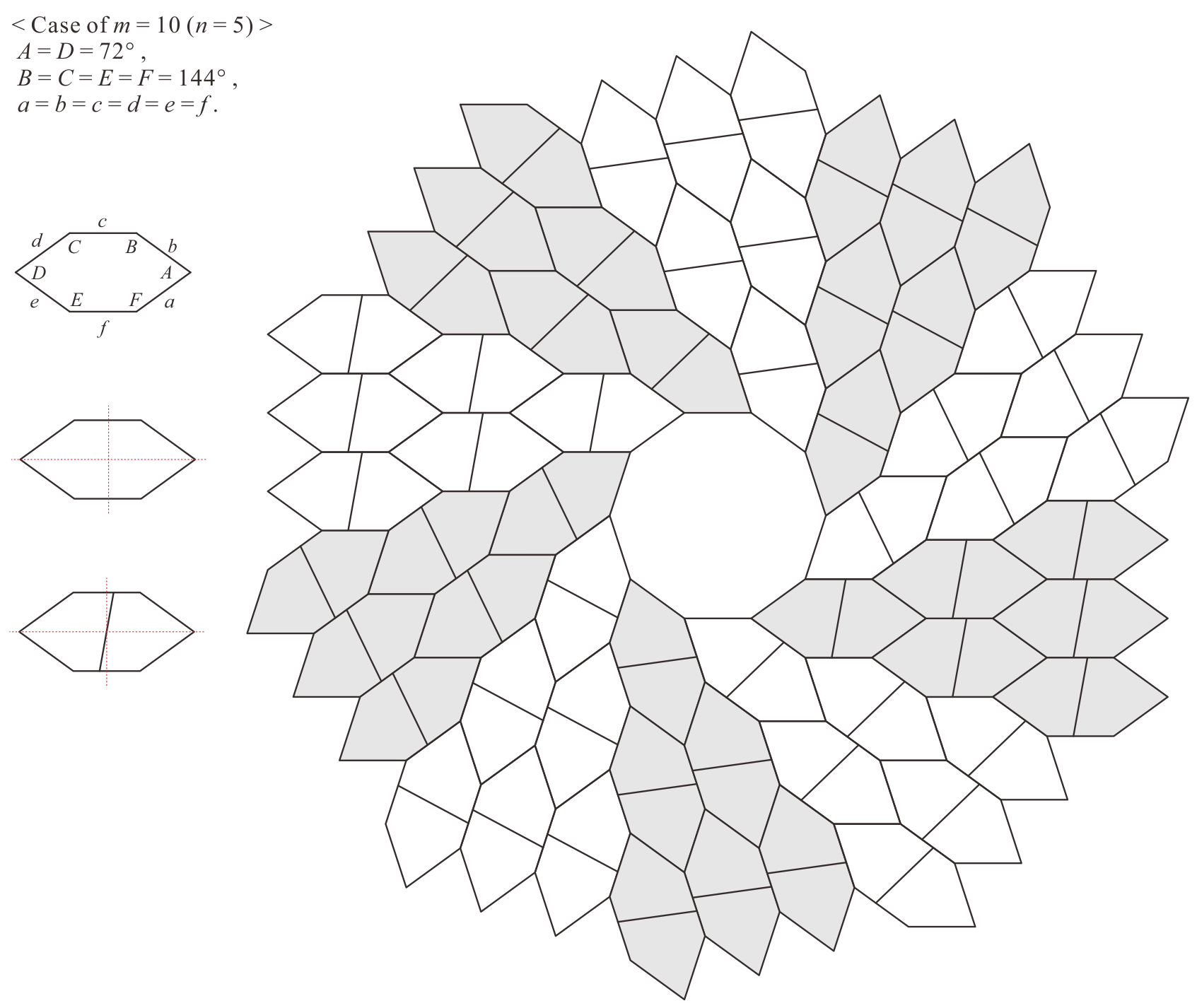} 
  \caption{{\small 
Rotationally symmetric tiling with $C_{10}$ symmetry with a regular 
convex 10-gonal hole at the center by a convex pentagon
}
\label{fig25}
}
\end{figure}

\begin{figure}[htbp]
 \centering\includegraphics[width=11.5cm,clip]{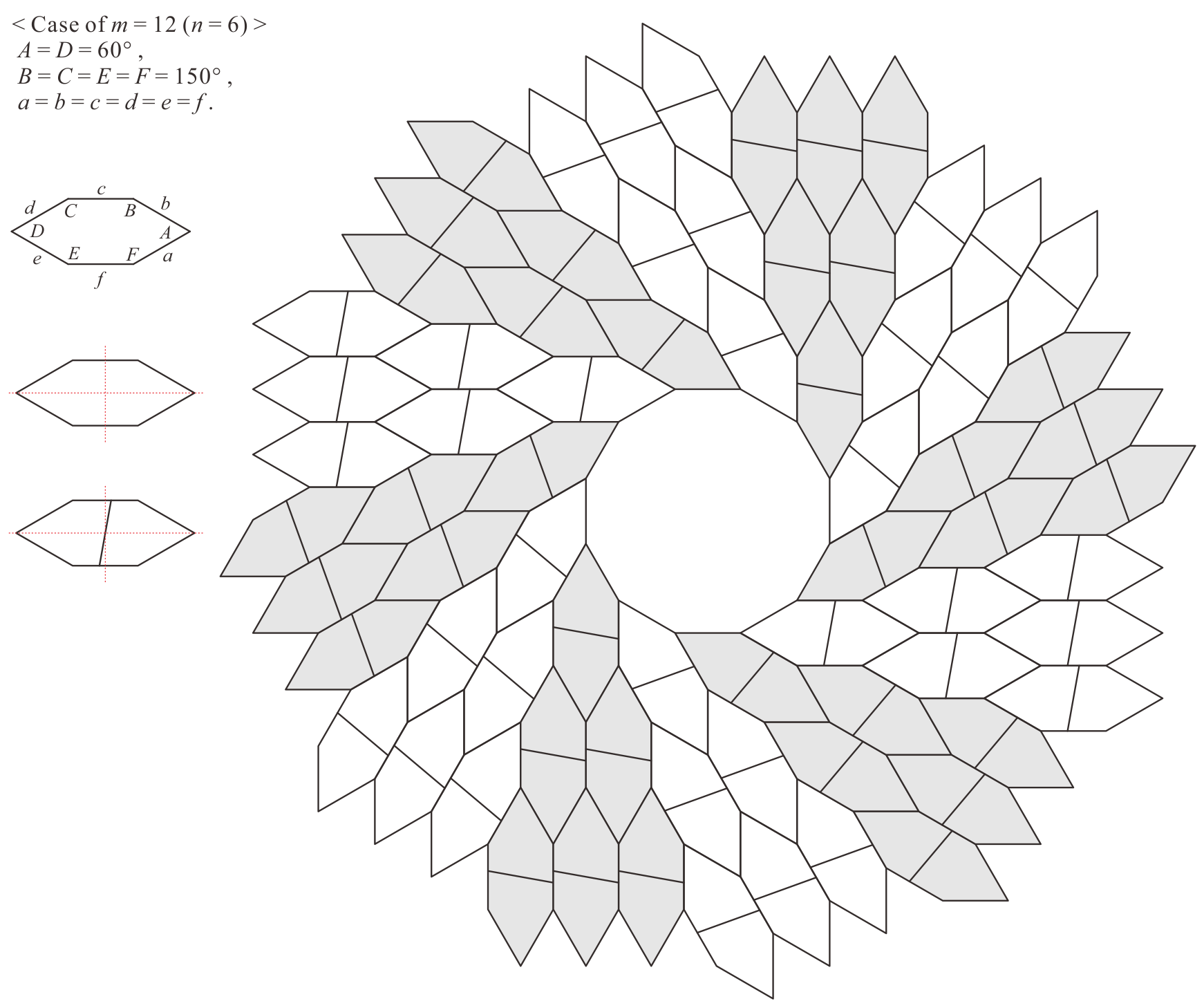} 
  \caption{{\small 
Rotationally symmetric tiling with $C_{12}$ symmetry with a regular 
convex 12-gonal hole at the center by a convex pentagon
}
\label{fig26}
}
\end{figure}

\begin{figure}[htbp]
 \centering\includegraphics[width=11.5cm,clip]{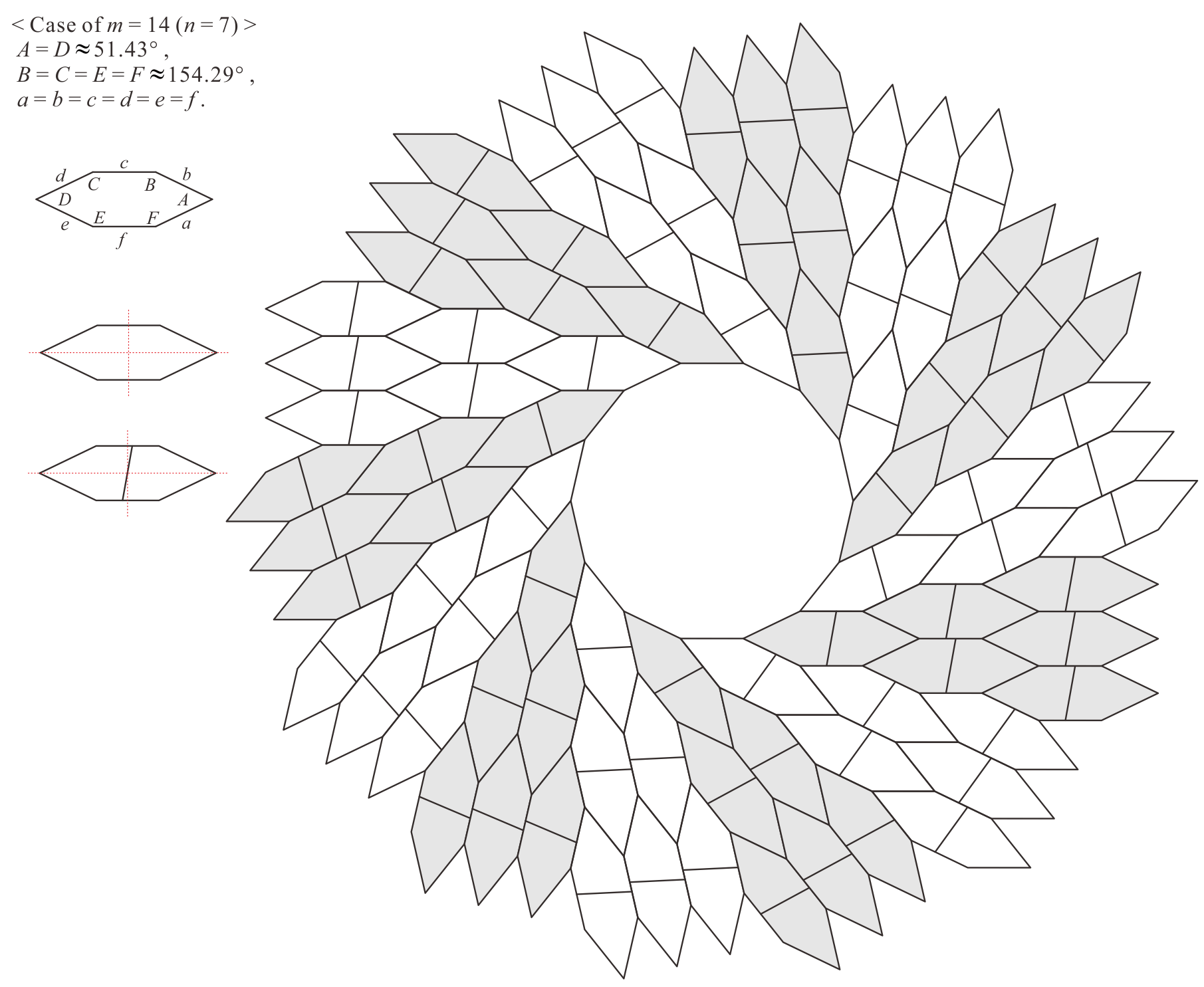} 
  \caption{{\small 
Rotationally symmetric tiling with $C_{14}$ symmetry with a regular 
convex 14-gonal hole at the center by a convex pentagon
}
\label{fig27}
}
\end{figure}

\begin{figure}[htbp]
 \centering\includegraphics[width=11.5cm,clip]{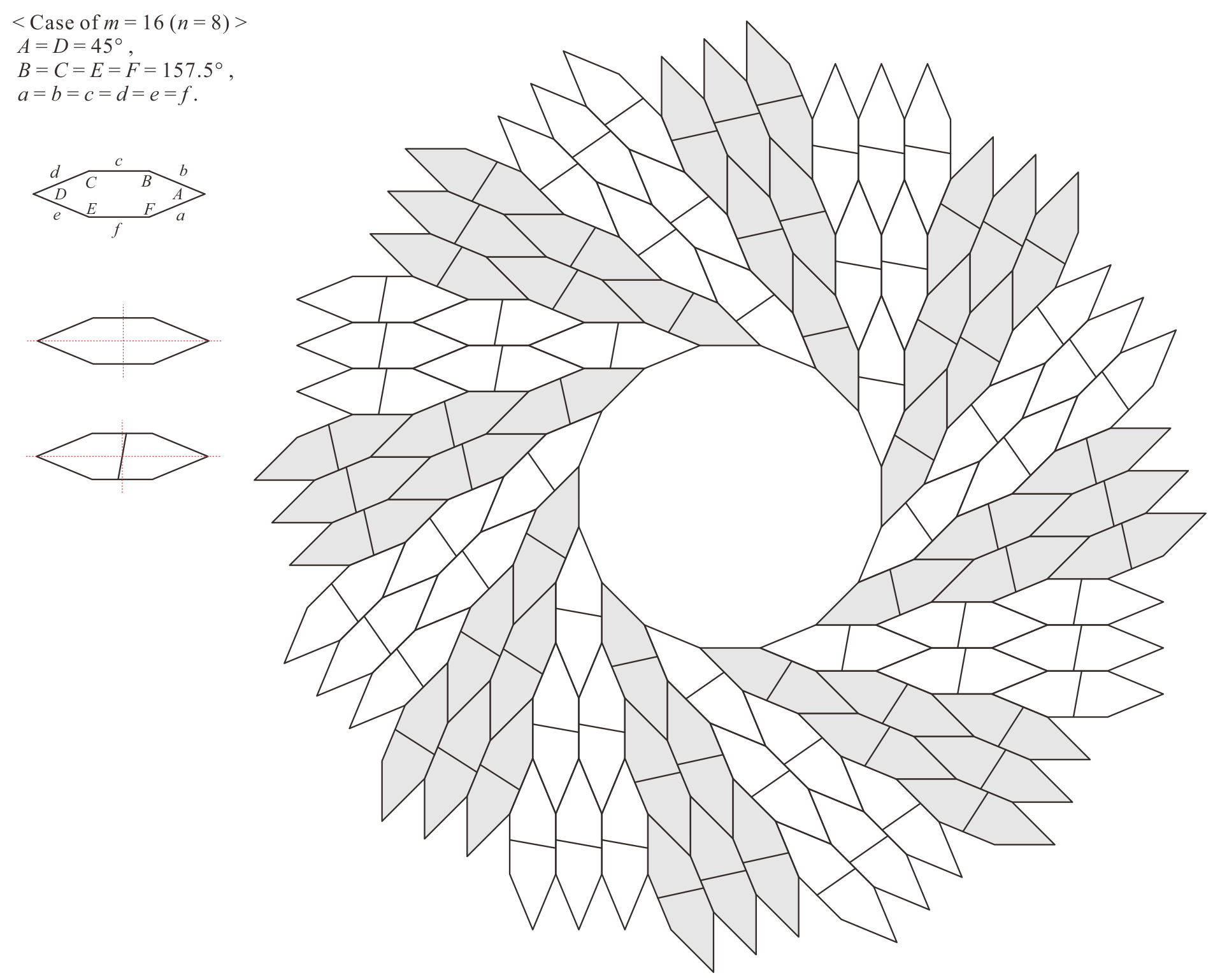} 
  \caption{{\small 
Rotationally symmetric tiling with $C_{16}$ symmetry with a regular 
convex 16-gonal hole at the center by a convex pentagon
}
\label{fig28}
}
\end{figure}

Because the convex hexagonal tile that satisfies ``$A = D = 90^ \circ ,\;
B = E = 120^ \circ ,\;C = F = 150^ \circ ,\; a = b = c = d = e = f$ " shown in 
Figures~\ref{fig08} and \ref{fig12} is the equilateral convex parallelohexagon 
with $\frac{360^ \circ }{3} = 120^ \circ $ and $\frac{360^ \circ }{4} = 90^ \circ $, 
it can form rotationally symmetric tilings with $C_{3}$ or $C_{4}$ symmetry, a 
rotationally symmetric tiling with $C_{3}$ symmetry with an equilateral convex 
hexagonal hole with $D_{3}$ symmetry, and a rotationally symmetric tiling with 
$C_{4}$ symmetry with an equilateral convex octagonal hole with $D_{4}$ 
symmetry at the center.  As an equilateral convex parallelohexagons 
with $C_{2}$ symmetry that has the same properties as above, there is a 
convex hexagonal tile that satisfies ``$A = D = 72^ \circ,\;B = E = 120^ \circ,\;
C = F = 168^ \circ,\;a = b = c = d = e = f$." This convex hexagonal tile can 
form rotationally symmetric tilings with $C_{5}$ or $C_{3}$ symmetry 
(see Figures~\ref{fig14}(a) and \ref{fig14}(b)), a rotationally symmetric tiling 
with $C_{5}$ symmetry with an equilateral convex 10-gonal hole with 
$D_{5}$ symmetry (see Figure~\ref{fig14}(c)),  and a rotationally symmetric 
tiling with $C_{3}$ symmetry with an equilateral convex hexagonal hole with 
$D_{3}$ symmetry at the center (see Figure~\ref{fig14}(d)). 
Note that Figure~\ref{fig14} shows the tilings with convex pentagons based on 
the convex hexagons.

\renewcommand{\figurename}{{\small Figure.}}

The convex hexagonal tiles that satisfy (\ref{eq4}) and have particular angles 
(i.e., angles corresponding to $\frac{360^ \circ }{n})$ of two or more 
types can form multiple rotationally symmetric tilings with an equilateral 
polygonal hole at the center, and multiple rotationally symmetric tilings. 
Such equilateral convex parallelohexagons with $C_{2}$ symmetry that can form 
multiple rotationally symmetric tilings with an equilateral polygonal hole 
at the center and multiple rotationally symmetric tilings are the above two 
cases with angles of ``$90^ \circ ,\;120^ \circ ,\;150^ \circ $" and 
``$72^ \circ ,\;120^ \circ ,\;168^ \circ $." It is because that interior angles of 
the convex hexagons that satisfy (\ref{eq4}) can be selected up to two different 
values corresponding to $\frac{360^ \circ }{n}$ when $n$ is an integer 
greater than or equal to three. Therefore, the two particular angles of this 
special case can be selected from ``$72^ \circ ,\;90^ \circ ,\;120^ \circ $."

\renewcommand{\figurename}{{\small Figure.}}

%%%%%%%%%%%%%%%%%%%%%%%%%%%%%%%%%%%%%%%%%%%%%%%%%%%%%%%%%%%%%%%%%%%%%%
%%%%%%%%%%%%%%%%%%%%%%%%%%%%%%%%%%%%%%%%%%%%%%%%%%%%%%%%%%%%%%%%%%%%%%

\section{Conclusions}
\label{section5}

This study summarizes how the convex pentagons and hexagons can form 
rotationally symmetric tilings based on the information in \cite{Iliev_2018} 
and \cite{Klaassen_2016}. The $n$-fold rotationally symmetric 
edge-to-edge tilings formed of convex hexagonal tiles that satisfy (\ref{eq1}) 
can be divided into $C_{n}$ symmetry and $D_{n}$ symmetry depending on 
the selection of the interior angles of the vertices in the convex hexagon. The convex 
pentagonal tiles generated by bisecting the convex hexagonal tiles can form 
rotationally symmetric non-edge-to-edge tilings with $C_{n}$ or $D_{n}$ symmetry 
depending on the shape or division method of the convex hexagon. In contrast, 
the $n$-fold rotationally symmetric edge-to-edge tilings with the convex 
pentagonal tiles shown in \cite{Sugimoto_2020_1} can form only $C_{n}$ symmetry. 
In addition, this study demonstrated a convex hexagons (convex pentagons) that 
can form rotationally symmetric tilings with an equilateral convex polygonal hole 
at the center.

In \cite{Iliev_2018}, there are cases where $m = 6, 8, 12$ present in Table~\ref{tab2} 
and tilings as shown in Figure~\ref{fig08} are described. However, there is no description 
of properties (tiles conditions such as (\ref{eq1}), (\ref{eq2}), and (\ref{eq3})) 
as discussed in this study. In \cite{Klaassen_2016}, there is no description 
of tilings with an equilateral convex polygonal hole at the center, and there 
is insufficient description regarding the convex hexagonal tile (such as 
the tile condition).

%%%%%%%%%%%%%%%%%%%%%%%%%%%%%%%%%%%%%%%%%%%%%%%%%%%%%%%%%%%%%%%%%%%%%%
%%%%%%%%%%%%%%%%%%%%%%%%%%%%%%%%%%%%%%%%%%%%%%%%%%%%%%%%%%%%%%%%%%%%%%

%%%%%%%%%%%%%%%%%%%%%%%%%%%%%%%%%%%%%%%%%%%%%%%%%%%%%%%%%%%%%%%%%%%%%%
%%%%%%%%%%%%%%%%%%%%%%%%%%%%%%%%%%%%%%%%%%%%%%%%%%%%%%%%%%%%%%%%%%%%%%

\renewcommand{\figurename}{{\small Figure.}}

\renewcommand{\figurename}{{\small Figure.}}

\renewcommand{\figurename}{{\small Figure.}}

\renewcommand{\figurename}{{\small Figure.}}

\renewcommand{\figurename}{{\small Figure.}}

\renewcommand{\figurename}{{\small Figure.}}

\renewcommand{\figurename}{{\small Figure.}}

\renewcommand{\figurename}{{\small Figure.}}

\renewcommand{\figurename}{{\small Figure.}}

\renewcommand{\figurename}{{\small Figure.}}

\renewcommand{\figurename}{{\small Figure.}}

\renewcommand{\figurename}{{\small Figure.}}

\renewcommand{\figurename}{{\small Figure.}}

\renewcommand{\figurename}{{\small Figure.}}

\end{document}